\documentclass{article}
\usepackage{ntheorem}
\usepackage{amsmath}
\usepackage{amssymb}
\usepackage{amstext}
\usepackage{xspace}
\usepackage{xy}
\xyoption{all}
\usepackage{hyperref}
\newcommand{\pullback}[1][dr]{\save*!/#1-1.2pc/#1:(-1,1)@^{|-}\restore}
\makeatother
 \newdir{ >}{{}*!/-7.5pt/@{>}}
 \newdir{|>}{!/4.5pt/@{|}*:(1,-.2)@^{>}*:(1,+.2)@_{>}}
 \newdir{ |>}{{}*!/-3pt/@{|}*!/-7.5pt/:(1,-.2)@^{>}*!/-7.5pt/:(1,+.2)@_{>}}
\newcommand{\xyline}[2][]{\ensuremath{\smash{\xymatrix@1#1{#2}}}}
\newcommand{\xyinline}[2][]{\ensuremath{\smash{\xymatrix@1#1{#2}}}^{\rule[8.5pt]{0pt}{0pt}}}
\makeatletter

\newif\ifignore % when set to true, additional text appears containing
\ignorefalse
\newcommand{\auxproof}[1]{
\ifignore\mbox{}\newline
\textbf{PROOF:} \dotfill\newline
{\it #1}\mbox{}\newline
\textbf{ENDPROOF}\dotfill
\fi}

\newcommand{\QEDbox}{\square}
\newcommand{\QED}{\hspace*{\fill}$\QEDbox$}

\newcommand{\after}{\mathrel{\circ}}
\newcommand{\Cat}[1]{\ensuremath{\mathbf{#1}}}
\newcommand{\field}[1]{\ensuremath{\mathbb{#1}}}
\newcommand{\op}{\ensuremath{^{\mathrm{op}}}}
\newcommand{\idmap}[1][]{\ensuremath{\mathrm{id}_{#1}}}
\renewcommand{\Im}{\ensuremath{\mathrm{Im}}}
\newcommand{\coker}{\ensuremath{\mathrm{coker}}}
\newcommand{\charac}{\ensuremath{\mathrm{char}}}

\newcommand{\Dom}{\ensuremath{\mathrm{Dom}}}

\newcommand{\inprod}[2]{\ensuremath{\langle #1\,|\,#2 \rangle}}

\newcommand{\Sub}{\ensuremath{\mathrm{Sub}}}
\newcommand{\KSub}{\ensuremath{\mathrm{KSub}}}
\newcommand{\orthogonal}{\mathrel{\bot}}

\newcommand{\sasaki}{\mathrel{\supset}}
\newcommand{\andthen}{\mathrel{\&}}
\newcommand{\tuple}[2]{\ensuremath{\protect\langle #1\,,\,#2 \protect\rangle}}

\newcommand{\Rel}{\Cat{Rel}\xspace}
\newcommand{\PInj}{\Cat{PInj}\xspace}
\newcommand{\Hilb}{\Cat{Hilb}\xspace}
\newcommand{\PHilb}{\Cat{PHilb}\xspace}
\newcommand{\Sets}{\Cat{Sets}\xspace}
\newcommand{\DCK}{\Cat{DagKerCat}\xspace}

\newcommand{\set}[2]{\{#1\;|\;#2\}}
\newcommand{\setin}[3]{\{#1\in#2\;|\;#3\}}
\newcommand{\conjun}{\mathrel{\wedge}}
\newcommand{\disjun}{\mathrel{\vee}}
\newcommand{\all}[2]{\forall_{#1}.\,#2}
\newcommand{\allin}[3]{\forall_{#1\in#2}.\,#3}
\newcommand{\ex}[2]{\exists_{#1}.\,#2}
\newcommand{\exin}[3]{\exists_{#1\in#2}.\,#3}
\newcommand{\fibno}[2]{\raisebox{.00in}
           {\mbox{ ${{\raisebox{-.05in}{$\scriptstyle #1$}\atop
           {\scriptscriptstyle\downarrow}}
           \atop{\scriptstyle #2}}$}}}
\newcommand{\downset}{\mathop{\downarrow}}
\newcommand{\effect}[1]{\mathfrak{E}(#1)}

\theoremnumbering{arabic}
\theoremstyle{plain}
\theorembodyfont{\itshape}
\theoremheaderfont{\normalfont\bfseries}
\theoremseparator{}
\newtheorem{theorem}{Theorem}
\newtheorem{lemma}[theorem]{Lemma}
\newtheorem{proposition}[theorem]{Proposition}
\newtheorem{corollary}[theorem]{Corollary}
\theorembodyfont{\normalfont}
\newtheorem{definition}[theorem]{Definition}
\newtheorem{example}[theorem]{Example}

\theoremstyle{nonumberplain}
\theoremheaderfont{\scshape}

\newtheorem{proof}{Proof}

\qedsymbol{\ensuremath{\Box}}

\begin{document}

\title{Quantum Logic in Dagger Kernel Categories}
\author{Chris Heunen \and Bart Jacobs}

\maketitle

\begin{abstract}
This paper investigates quantum logic from the perspective of
categorical logic, and starts from minimal assumptions, namely the
existence of involutions/daggers and kernels. The resulting structures
turn out to (1)~encompass many examples of interest, such as
categories of relations, partial injections, Hilbert spaces (also
modulo phase), and Boolean algebras, and (2)~have interesting
categorical/logical/order-theoretic properties, in terms of kernel
fibrations, such as existence of pullbacks, factorisation,
orthomodularity, atomicity and completeness. For instance, the Sasaki
hook and and-then connectives are obtained, as adjoints, via the
existential-pullback adjunction between fibres.
\end{abstract}
%\end{frontmatter}

\section{Introduction}\label{IntroSec}

Dagger categories \Cat{D} come equipped with a special functor
$\dag\colon \Cat{D}\op \rightarrow \Cat{D}$ with $X^\dag=X$ on objects
and $f^{\dag\dag}=f$ on morphisms. A simple example is the category
\Rel of sets and relations, where $\dag$ is reversal of relations. A
less trivial example is the category \Hilb of Hilbert spaces and
continuous linear transformations, where $\dag$ is induced by the
inner product. The use of daggers, mostly with additional assumptions,
dates back to~\cite{MacLane61,Puppe62}. Daggers are currently of
interest in the context of quantum
computation~\cite{AbramskyC04,Selinger07,CoeckeP06}. The dagger
abstractly captures the reversal of a computation.

Mostly, dagger categories are used with fairly strong additional
assumptions, like compact closure in~\cite{AbramskyC04}. Here we wish
to follow a different approach and start from minimal
assumptions. This paper is a first step to understand quantum logic,
from the perspective of categorical logic (see
\textit{e.g.}~\cite{MakkaiR77,KockR77,Taylor99,Jacobs99a}).  It grew
from the work of one of the authors~\cite{Heunen08b}.  Although that
paper enjoys a satisfactory relation to traditional quantum
logic~\cite{Heunen09b}, this one generalises it, by taking the notion
of dagger category as starting point, and adding kernels, to be used
as predicates. The interesting thing is that in the presence of a
dagger $\dag$ much else can be derived.  As usual, it is quite subtle
what precisely to take as primitive. A referee identified the
reference~\cite{Crown75} as an earlier precursor to this work. It
contains some crucial ingredients, like orthomodular posets of dagger
kernels, but without the general perspective given by categorical
logic.

%% The paper basically contains two parts, which produce roughly
%% the following two results.
%% \begin{itemize}
%% \item Kernels in a dagger category yield a factorisation system
%% of ``zero-epis'' and kernels, satisfying diagonal fill-in.

%% \item Assuming a suitable tensor sum $\oplus$ yields a (partial)
%% sum $\amalg$ and a partial order on homsets.
%% \end{itemize}

Upon this structure of ``dagger kernel categories'' the paper
constructs pullbacks of kernels and factorisation (both similar
to~\cite{Freyd64}). It thus turns out that the kernels form a
``bifibration'' (both a fibration and an opfibration,
see~\cite{Jacobs99a}). This structure can be used as a basis for
categorical logic, which captures substitution in predicates by
reindexing (pullback) $f^{-1}$ and existential quantification by
op-reindexing $\exists_f$, in such a way that $\exists_{f} \dashv
f^{-1}$. From time to time we use fibred terminology in this paper,
but familiarity with this fibred setting is not essential. We find
that the posets of kernels (fibres) are automatically orthomodular
lattices~\cite{Kalmbach83}, and that the Sasaki hook and and-then
connectives appear naturally from the existential-pullback
adjunction. Additionally, a notion of Booleanness is identified for
these dagger kernel categories. It gives rise to a generic
construction that generalises how the category of partial injections
can be obtained from the category of relations.

Apart from this general theory, the paper brings several important
examples within the same setting---of dagger kernel
categories. Examples are the categories \Rel and \PInj of relations
and partial injections. Additionally, the category \Hilb is an
example---and, interestingly---also the category \PHilb of Hilbert
spaces modulo phase. The latter category provides the framework in
which physicists typically work~\cite{Lahti04}. It has much weaker
categorical structure than \Hilb. We also present a construction to
turn an arbitrary Boolean algebra into a dagger kernel category.

The authors are acutely aware of the fact that several of the example
categories have much richer structure, involving for instance a tensor
sum $\oplus$ and a tensor product $\otimes$ with associated scalars
and traced monoidal structure. This paper deliberately concentrates
solely on (the logic of) kernels. There are interesting differences
between our main examples: for instance, \Rel and \PInj are Boolean,
but \Hilb is not; in \PInj and \Hilb ``zero-epis'' are epis, but not
in \Rel; \Rel and \Hilb have biproducts, but \PInj does not.

The paper is organised as follows. After introducing the notion of
dagger kernel category in Section~\ref{KernelSec}, the main examples
are described in Section~\ref{ExamplesSec}. Factorisation and
(co)images occur in Sections~\ref{FactorisationSec}
and~\ref{ImageCoimageSec}. Section~\ref{LogicSec} introduces the
Sasaki hook and and-then connectives via adjunctions, and investigates
Booleanness. Finally, Sections~\ref{HomsetOrderSec}
and~\ref{ComplAtomSec} investigate some order-theoretic aspects of
homsets and of kernel posets (atomicity and completeness).

A follow-up paper~\cite{Jacobs09a} introduces a new category
\Cat{OMLatGal} of orthomodular lattices with Galois connections
between them, shows that it is a dagger kernel category, and that
every dagger kernel category \Cat{D} maps into it via the kernel
functor $\KSub\colon\Cat{D} \rightarrow \Cat{OMLatGal}$, preserving
the dagger kernel structure. This gives a wider context.

\section{Daggers and kernels}\label{KernelSec}

Let us start by introducing the object of study of this paper.

\begin{definition}
\label{DagCatKerDef}
A \emph{dagger kernel category} consists of:
\begin{enumerate}
\item a dagger category \Cat{D}, with dagger $\dag\colon \Cat{D}\op
  \rightarrow \Cat{D}$;

\item a zero object 0 in \Cat{D};

\item kernels $\ker(f)$ of arbitrary maps $f$ in \Cat{D}, which are 
dagger monos.
\end{enumerate}
A \emph{morphism of dagger kernel categories} is a functor $F$
preserving the relevant structure:
\begin{enumerate}
\item $F(f^\dag) = F(f)^\dag$;
\item $F(0)$ is again a zero object;
\item $F(k)$ is a kernel of $F(f)$ if $k$ is a kernel of $f$.
\end{enumerate}
Dagger kernel categories and their morphisms form a category $\DCK$.
\end{definition}

\begin{definition}
A dagger kernel category is called Boolean if $m\conjun n 
=  0$ implies $m^{\dag} \after n = 0$, for all kernels $m,n$. 
\end{definition}

The name Boolean
will be explained in Theorem~\ref{BooleanLem}. We shall
later rephrase the Booleanness condition as: kernels are disjoint if
and only if they are orthogonal, see Lemma~\ref{KerLem}. \\

The dagger $\dag$ satisfies $X^{\dag} = X$ on objects and
$f^{\dag\dag} = f$ on morphisms. It comes with a number of
definitions. A map $f$ in \Cat{D} is called a dagger mono(morphism) if
$f^{\dag} \after f = \idmap$ and a dagger epi(morphism) if $f\after
f^{\dag} = \idmap$. Hence $f$ is a dagger mono if and only if
$f^{\dag}$ is a dagger epi. A map $f$ is a dagger iso(morphism) when
it is both dagger monic and dagger epic; in that case $f^{-1} =
f^{\dag}$ and $f$ is sometimes called unitary (in analogy with Hilbert
spaces). An endomorphism $p\colon X\rightarrow X$ is called
self-adjoint if $p^{\dag} = p$.

%% A projection $p$ is a self-adjoint idempotent: $p^{\dag} = p = p
%% \after p$.

The zero object $0\in\Cat{D}$ is by definition both initial and
final. Actually, in the presence of $\dag$, initiality implies
finality, and vice-versa. For an arbitrary object $X\in\Cat{D}$, the
unique map $X\rightarrow 0$ is then a dagger epi and the unique map
$0\rightarrow X$ is a dagger mono. The ``zero'' map $0 = 0_{X,Y} =
(X\rightarrow 0 \rightarrow Y)$ satisfies $(0_{X,Y})^{\dag} =
0_{Y,X}$. Notice that $f \after 0 = 0 = 0 \after g$. Usually there is
no confusion between 0 as zero object and 0 as zero map.  Two maps
$f\colon X\rightarrow Z$ and $g\colon Y\rightarrow Z$ with common
codomain are called orthogonal, written as $f\orthogonal g$, if
$g^{\dag} \after f = 0$---or, equivalently, $f^{\dag} \after g = 0$.

Let us recall that a kernel of a morphism $f\colon X\rightarrow Y$ is
a universal morphism $k\colon\ker(f)\rightarrow X$ with $f\after k =
0$.  Universality means that for an arbitrary $g\colon Z\rightarrow X$
with $f\after g = 0$ there is a unique map $g'\colon Z\rightarrow
\ker(f)$ with $k \after g' = g$. Kernels are automatically (ordinary)
monos. Just like we write $0$ both for a zero object and for a zero
map, we often write $\ker(f)$ to denote either a kernel map, or the
domain object of a kernel map.

Definition~\ref{DagCatKerDef} requires that kernels are dagger
monos. This requirement involves a subtlety: kernels are closed under
arbitrary isomorphisms but dagger monos are only closed under dagger
isomorphisms. Hence we should be more careful in this
requirement. What we really mean in Definition~\ref{DagCatKerDef} is
that for each map $f$, among all its isomorphic kernel maps, there is
at least one dagger mono. We typically choose this dagger mono as
representant $\ker(f)$ of the equivalence class of kernel maps.

We shall write $\KSub(X)$ for the poset of (equivalence classes) of
kernels with codomain $X$. The order $(M \rightarrowtail X) \leq
(N\rightarrowtail X)$ in $\KSub(X)$ is given by the presence of a
(necessarily unique) map $M\rightarrow N$ making the obvious triangle
commute.  Intersections in posets like $\KSub(X)$, if they exist, are
given by pullbacks, as in:
$$\xymatrix{
\bullet\ar@{ >->}[r]\ar@{ >->}[d]\ar@{ >->}[dr]|{m\conjun n} & 
   M\ar@{ >->}[d]^{m} \\
N\ar@{ >->}[r]_-{n} & X\rlap{.}
}$$

In presence of the dagger $\dag$, cokernels come for free: one can
define a cokernel $\coker(f)$ as $\ker(f^{\dag})^{\dag}$. Notice that
we now write $\ker(f)$ and $\coker(f)$ as morphisms.  This cokernel
$\coker(f)$ is a dagger epi. Finally, we define $m^{\perp} =
\ker(m^{\dag})$, which we often write as $m^\perp : M^\perp
\rightarrowtail X$ if $m\colon M \rightarrowtail X$. This notation is
especially used when $m$ is a mono. In diagrams we typically write a
kernel as $\smash{\xymatrix@C-.5pc{\ar@{|>->}[r] & }}$ and a cokernel
as $\xymatrix@C-.5pc{\ar@{-|>}[r] & }$.

The following Lemma gives some basic observations.

\begin{lemma}
\label{KerLem}
In a dagger kernel category,
\begin{enumerate}
\item $\ker(\smash{X\stackrel{0}{\rightarrow} Y}) =
  (\smash{X\stackrel{\idmap}{\rightarrow}X})$ and
  $\ker(\smash{X\stackrel{\idmap}{\rightarrow}X}) =
  (\smash{0\stackrel{0}{\rightarrow}X})$; these yield greatest and
  least elements $1,0\in\KSub(X)$, respectively;

\item $\ker(\ker(f)) = 0$;

\item $\ker(\coker(\ker(f))) = \ker(f)$, as subobjects;

\item $m^{\perp\perp} = m$ if $m$ is a kernel;

\item A map $f$ factors through $g^{\perp}$ iff $f\orthogonal g$ iff
  $g \orthogonal f$ iff $g$ factors through $f^{\perp}$; in particular
  $m \leq n^{\perp}$ iff $n\leq m^{\perp}$, for monos $m,n$; hence
  $(-)^{\perp} \colon \KSub(X)$
  $\smash{\stackrel{\cong}{\longrightarrow}}$ $\KSub(X)\op$;

\item if $m\leq n$, for monos $m,n$, say via $m = n\after \varphi$, then:
\begin{enumerate}
\item if $m,n$ are dagger monic, then so is $\varphi$;
\item if $m$ is a kernel, then so is $\varphi$.
\end{enumerate}

\item Booleanness amounts to $m\conjun n = 0 \Leftrightarrow
  m\orthogonal n$, \textit{i.e.}~disjointness is orthogonality, for
  kernels.
\end{enumerate}
\end{lemma}

\begin{proof}
We skip the first two points because they are obvious and start with the
third one. Consider the following diagram for an arbitrary $f\colon
X\rightarrow Y$:
$$\xymatrix{
\ker(f)\ar@{ |>->}[rr]^-{k}\ar@{..>}[d]<-1ex>_-{k'} & & 
   X\ar[rr]^-{f}\ar@{-|>}[drr]_-{c} & & Y \\
\ker(\coker(\ker(f)))\ar@{ |>->}[urr]_-{\ell}\ar@{..>}[u]<-1ex>_-{\ell'} 
   & & & & \coker(\ker(f)). \ar@{..>}[u]_{f'}
}$$

\noindent By construction $f\after k = 0$ and $c\after k=0$. Hence
there are $f'$ and $k'$ as indicated. Since $f \after \ell = f' \after
c \after \ell = f' \after 0 = 0$ one gets $\ell'$. Hence the kernels
$\ell$ and $k$ represent the same subobject.

For the fourth point, notice that if $m=\ker(f)$, then
$$m^{\perp\perp}
=
\ker(\ker(m^{\dag})^{\dag})
=
\ker(\coker(\ker(f)))
=
\ker(f)
=
m.$$

\noindent Next,
$$\begin{array}{rcl}
\mbox{$f$ factors through $g^{\perp}$}
& \Longleftrightarrow &
g^{\dag} \after f = 0 \\
& \Longleftrightarrow &
f^{\dag} \after g = 0
\hspace*{\arraycolsep}\Longleftrightarrow\hspace*{\arraycolsep}
\mbox{$g$ factors through $f^{\perp}$.}
\end{array}$$

\noindent If, in the sixth point, $m = n \after \varphi$ and $m,n$ are
dagger monos, then $\varphi^{\dag} \after \varphi = (n^{\dag} \after
m)^{\dag} \after \varphi = m^{\dag} \after n \after \varphi = m^{\dag}
\after m = \idmap$. And if $m = \ker(f)$, then $\varphi = \ker(f\after
n)$, since: (1)~$f \after n \after \varphi = f \after m = 0$, and
(2)~if $f\after n \after g = 0$, then there is a $\psi$ with $m \after
\psi = n \after g$, and this gives a unique $\psi$ with $\varphi
\after \psi = g$, where uniqueness of this $\psi$ comes from $\varphi$
being monic.

Finally, Booleanness means that $m\conjun n = 0$ implies $m^{\dag}
\after n = 0$, which is equivalent to $n^{\dag} \after m = 0$, which
is $m\orthogonal n$ by definition. The reverse implication is easy,
using that the meet $\conjun$ of monos is given by pullback: if
$m\after f = n \after g$, then $f = m^{\dag} \after m \after f =
m^{\dag} \after n \after g = 0 \after g = 0$. Similarly, $g=0$. Hence
the zero object $0$ is the pullback $m\conjun n$ of $m,n$. 
\end{proof}

Certain constructions from the theory of Abelian
categories~\cite{Freyd64} also work in the current setting. This
applies to the pullback construction in the next result, but also, to
a certain extent, to the factorisation of
Section~\ref{FactorisationSec}.

\begin{lemma}
\label{PullbackLem}
Pullbacks of kernels exist, and are kernels again. Explicitly,
given a kernel $n$ and map $f$ one obtains a pullback:
$$\begin{array}{rcrcl}
\raisebox{1.5em}{$\xymatrix{
M\ar[rr]^-{f'}\ar@{ |>->}[d]_{f^{-1}(n)}\pullback & & N\ar@{ |>->}[d]^{n} \\
X\ar[rr]_-{f} & & Y
}$}
& \quad\mbox{as}\quad &
f^{-1}(n)
& = &
\ker(\coker(n) \after f).
\end{array}$$

\noindent If $f$ is a dagger epi, so is $f'$.
\end{lemma}

By duality there are of course similar results about pushouts
of cokernels.

\begin{proof}
For convenience write $m$ for the dagger kernel $f^{-1}(n) =
\ker(\coker(n) \after f)$. By construction, $\coker(n) \after f \after
m = 0$, so that $f\after m$ factors through $\ker(\coker(n)) = n$, say
via $f'\colon M\rightarrow N$ with $n\after f' = f \after m$, as in
the diagram. This yields a pullback: if $a\colon Z\rightarrow X$ and
$b\colon Z\rightarrow N$ satisfy $f\after a = n\after b$, then
$\coker(n) \after f \after a = \coker(n) \after n \after b = 0 \after
b = 0$, so that there is a unique map $c\colon Z\rightarrow M$ with
$m\after c = a$. Then $f'\after c = b$ because $n$ is monic.

In case $f$ is dagger epic, $f \after f^\dag \after n = n$. Hence
there is a morphism $f''$ making following diagram commute, as the
right square is a pullback: 
\[\xymatrix@R-3ex{
  N \ar@{-->}_-{f''}[dr] \ar@(d,l)_-{f^\dag \after n}[dddr]
  \ar@(r,ul)^-{\idmap}[drrr] \\ 
  & M \pullback \ar^-{f'}[rr] \ar@{ |>->}_-{m}[dd] 
  && N \ar@{ |>->}^-{n}[dd] \\ \\
  & X \ar_-{f}[rr] && Y.  
}\]
Then $f'' = m^{\dag} \after m \after f'' = m^{\dag} \after f^{\dag} \after
n = (f \after m)^{\dag} \after n = (n \after f')^{\dag} \after n =
f'^{\dag} \after n^{\dag} \after n = f'^{\dag}$. Hence $f'$ is
dagger epic, too. 
% In case $f$ is a dagger epi we have $f \after f^{\dag} = \idmap$.
% Hence there are two adjacent pullbacks: \marginpar{why??}
% $$\xymatrix{
% N\ar[rr]^-{f''}\ar@{ |>->}[d]_{n}\pullback & &
%    M\ar[rr]^-{f'}\ar@{ |>->}[d]|{m = f^{-1}(n)}\pullback & & 
%    N\ar@{ |>->}[d]^{n} \\
% Y\ar[rr]_-{f^\dag} & & X\ar[rr]_-{f} & & Y
% }$$

% \noindent Then $f' \after f'' = \idmap$ because $n$ is monic. Further,
% $f'' = m^{\dag} \after m \after f'' = m^{\dag} \after f^{\dag} \after
% n = f'^{\dag} \after n^{\dag} \after n = f'^{\dag}$. Hence $f'$ is
% dagger epi. 
\end{proof}

\begin{corollary}
\label{PullbackCor} 
%% ~\\
Given these pullbacks of kernels we observe the following.
\begin{enumerate}
\item The mapping $X\mapsto \KSub(X)$ yields an indexed category
  $\Cat{D}\op\rightarrow\Cat{PoSets}$, using that each map $f\colon
  X\rightarrow Y$ in \Cat{D} yields a pullback (or substitution)
  functor $f^{-1}\colon \KSub(Y) \rightarrow \KSub(X)$. By the
  ``pullback lemma'', see \textit{e.g.}~\cite[Lemma 5.10]{Awodey06}
  or~\cite[III, 4, Exc.~8]{MacLane71}, such functors $f^{-1}$ preserve
  the order on kernels, and also perserve all meets (given by
  pullbacks). This (posetal) indexed category $\KSub \colon
  \Cat{D}\op\rightarrow\Cat{PoSets}$ forms a setting in which one can
  develop categorical logic for dagger categories, see
  Subsection~\ref{FibrationSubsec}.

% \item if kernels $n,m \in \KSub(X)$ satisfy $m \leq n$, say via $m = n
%   \after \varphi$, then the following diagram is a pullback,
%   \[\xymatrix{
%     M \pullback \ar@{=}[r] \ar@{ |>->}_-{\varphi}[d] & M \ar@{
%     |>->}^-{m}[d] \\
%     N \ar@{ |>->}_-{n}[r] & X,
%   }\]
%   so that $\varphi$ is a kernel too;

\item The following diagram is a pullback,
$$\xymatrix@R-.5pc{
\ker(f)\ar[rr]\ar@{ |>->}[d]\pullback & & 0\ar@{ |>->}[d] \\
X\ar[rr]_-{f} & & Y
}$$

\noindent showing that, logically speaking, falsum---\textit{i.e.}~the
bottom element $0\in\KSub(Y)$---is in general not preserved under
substitution. Also, negation/orthocomplementation $(-)^{\perp}$ does
not commute with substitution, because $1 = 0^{\perp}$ and $f^{-1}(1)
= 1$. 
\end{enumerate}
\end{corollary}

Being able to take pullbacks of kernels has some important
consequences.

\begin{lemma}
\label{KernelCompLem}
Kernels are closed under composition---and hence cokernels are, too.
\end{lemma}

\begin{proof}
We shall prove the result for cokernels, because it uses pullback
results as we have just seen. So assume we have (composable) cokernels
$e, d$; we wish to show $e \after d = \coker(\ker(e\after d))$. We
first notice, using Lemma~\ref{PullbackLem},
$$\ker(e \after d)
=
\ker(\coker(\ker(e)) \after d)
=
d^{-1}(\ker(e)),$$

\noindent yielding a pullback:
$$\xymatrix{
& & A\ar[rr]^-{d'}\ar@{ |>->}[d]|{m = \ker(e\after d)}\pullback & &
   B\ar@{ |>->}[d]^{\ker(e)} \\
K\ar@{ |>->}[rr]^-{\ker(d)}\ar@{-->}[urr]^{\varphi} & & 
   X\ar@{-|>}[rr]^-{d} & & 
   D\ar@{-|>}[rr]^-{e} & & E.
}$$

\noindent We intend to prove $e\after d = \coker(m)$. Clearly, $e
\after d \after m = e \after \ker(e) \after d' = 0 \after d' = 0$.
And if $f\colon X\rightarrow Y$ satisfies $f\after m = 0$, then $f
\after \ker(d) = f \after m \after \varphi = 0$, so because $d =
\coker(\ker(d))$ there is $f'\colon D\rightarrow Y$ with $f' \after d
= f$. But then: $f' \after \ker(e) \after d' = f' \after d \after m =
f\after m = 0$.  Then $f' \after \ker(e) = 0$, because $d'$ is
dagger epi because $d$ is, see Lemma~\ref{PullbackLem}. This
finally yields $f'' \colon E\rightarrow Y$ with $f'' \after e = f'$.
Hence $f'' \after e \after d = f$. 
\end{proof}

As a result, the logic of kernels has intersections, preserved by
substitution.  More precisely, the indexed category $\KSub(-)$ from
Corollary~\ref{PullbackCor} is actually a functor $\KSub \colon
\Cat{D}\op \rightarrow \Cat{MSL}$ to the category \Cat{MSL} of meet
semi-lattices. Each poset $\KSub(X)$ also has disjunctions, by
$m\disjun n = (m^{\perp} \conjun n^{\perp})^{\perp}$, but they are not
preserved under substitution/pullback $f^{-1}$. Nevertheless, $m
\disjun m^{\perp} = (m^{\perp} \conjun m^{\perp\perp})^{\perp} =
(m^{\perp} \conjun m)^{\perp} = 0^{\perp} = 1$.

The essence of the following result goes back to~\cite{Crown75}.

\begin{proposition}
\label{OrthmodularityProp}
Orthomodularity holds: for kernels $m\leq n$, say via $\varphi$
with $n\after \varphi = m$, one has pullbacks:
$$\xymatrix{
M\ar@{ |>->}[rr]^-{\varphi}\ar@{=}[d]\pullback & &
   N\ar@{ |>->}[d]^{n} & & 
   P\ar[d]\ar@{ |>->}[ll]_-{\varphi^{\perp}}\pullback[ld] \\
M\ar@{ |>->}[rr]_-{m} & & X & & M^{\perp}\ar@{ |>->}[ll]^-{m^{\perp}}
}$$

\noindent This means that $m \disjun (m^{\perp} \conjun n) = n$.
\end{proposition}

\begin{proof}
The square on the left is obviously a pullback.  For the one on the
right we use a simple calculation, following Lemma~\ref{PullbackLem}:
$$\begin{array}{rcll}
n^{-1}(m^{\perp})
& = &
\ker(\coker(m^{\perp}) \after n) \\
& = &
\ker(\coker(\ker(m^{\dag})) \after n) \\
& = &
\ker(m^{\dag} \after n)
   & \mbox{since $m^{\dag}$ is a cokernel} \\
& \smash{\stackrel{(*)}{=}} &
\ker(\varphi^{\dag}) \\
& = &
\varphi^{\perp},
\end{array}$$

\noindent where the marked equation holds because $n\after \varphi = m$,
so that $\varphi = n^{\dag} \after n \after \varphi = n^{\dag} \after m$
and thus $\varphi^{\dag} = m^{\dag} \after n$. Then:
$$m \disjun (m^{\perp} \conjun n)
\hspace*{\arraycolsep} = \hspace*{\arraycolsep}
(n \after \varphi) \disjun (n\after \varphi^{\perp}) \\
\hspace*{\arraycolsep}\smash{\stackrel{(*)}{=}}\hspace*{\arraycolsep}
n \after (\varphi \disjun \varphi^{\perp}) \\
\hspace*{\arraycolsep} = \hspace*{\arraycolsep}
n \after \idmap \\
\hspace*{\arraycolsep} = \hspace*{\arraycolsep}
n.$$

\noindent The (newly) marked equation holds because $n\after (-)$
preserves joins, since it is a left adjoint: $n\after k \leq m$ iff $k
\leq n^{-1}(m)$, for kernels $k,m$. 
\end{proof}

The following notion does not seem to have an established terminology,
and therefore we introduce our own.

\begin{definition}
\label{ZeroMonoEpiDef}
In a category with a zero object, a map $m$ is called a zero-mono if
$m\after f = 0$ implies $f=0$, for any map $f$. Dually, $e$ is
zero-epi if $f \after e = 0$ implies $f=0$.
In diagrams we write $\smash{\xymatrix@C-.5pc{\ar@{ >->}|{\circ}[r] & }}$
for zero-monos and $\smash{\xymatrix@C-.5pc{\ar@{->>}|{\circ}[r] & }}$ for
zero-epis.
\end{definition}

Clearly, a mono is zero-mono, since $m\after f = 0 = m\after 0$
implies $f=0$ if $m$ is monic. The following points are worth
making explicit.

\begin{lemma}
\label{ZeroMonoEpiLem}
In a dagger kernel category,
\begin{enumerate}
\item $m$ is a zero-mono iff $\ker(m) = 0$ and $e$ is a zero-epi iff
$\coker(e) = 0$;

\item $\ker(m\after f) = \ker(f)$ if $m$ is a zero-mono, and similarly,
$\coker(f \after e) = \coker(f)$ if $e$ is a zero-epi;

\item a kernel which is zero-epic is an isomorphism. \QED
\end{enumerate}
\end{lemma}

We shall mostly be interested in zero-epis (instead of zero-monos),
because they arise in the factorisation of
Section~\ref{FactorisationSec}. In the presence of dagger equalisers,
zero-epis are ordinary epis.  This applies to \Hilb and
\PInj. This fact is not really used, but is included because it
gives a better understanding of the situation. A \emph{dagger
equaliser category} is a dagger category that has equalisers which are
dagger monic.

\begin{lemma}
\label{EqualiserZeroEpiLem}
In a dagger equaliser category \Cat{D} where every dagger mono is a
kernel, zero-epis in \Cat{D} are ordinary epis. 
\end{lemma}
\begin{proof}
Assume a zero-epi $e\colon E\rightarrow X$ with two maps $f,g\colon
X\rightarrow Y$ satisfying $f\after e = g\after e$. We need to
prove $f=g$. Let $m\colon M\rightarrowtail X$ be the equaliser of 
$f,g$, with $h = \coker(m)$, as in:
$$\xymatrix{
E\ar@{->>}|{\circ}[rr]^-{e}\ar@{-->}[d]_{\varphi} & & 
   X \ar[rr]<.5ex>^-{f}\ar[rr]<-.5ex>_-{g} \ar@{-|>}[d]^(0.6){h=\coker(m)} & & Y \\
M\ar@{ |>->}[urr]_-{m} & & Z
}$$

\noindent This $e$ factors through the equaliser $m$, as indicated,
since $f\after e = g\after e$. Then $h \after e = h \after m \after
\varphi = 0 \after \varphi = 0$. Hence $h = 0$ because $e$ is
zero-epi. But $m$, being a dagger mono, is a dagger kernel. Hence $m =
\ker(\coker(m)) = \ker(h) = \ker(0) = \idmap$, so that
$f=g$. 
\end{proof}

\subsection{Indexed categories and fibrations}\label{FibrationSubsec}

The kernel posets $\KSub(X)$ capture the predicates on an object $X$,
considered as underlying type, in a dagger kernel category
\Cat{D}. Such posets are studied systematically in categorical logic,
often in terms of indexed categories
$\Cat{D}\op\rightarrow\Cat{Posets}$ or even as a so-called fibration
$\Big(\fibno{\KSub(\Cat{D})}{\Cat{D}}\Big)$, see~\cite{Jacobs99a}.  We
shall occasionally borrow terminology from this setting, but will not
make deep use of it. A construction that is definitely useful in the
present setting is the ``total'' category $\KSub(\Cat{D})$. It has
(equivalence classes of) kernels $M\rightarrowtail X$ as
objects. Morphisms $\smash{(M\stackrel{m}{\rightarrowtail} X)}
\longrightarrow \smash{(N\stackrel{n}{\rightarrowtail} Y)}$ in
$\KSub(\Cat{D})$ are maps $f\colon X \rightarrow Y$ in \Cat{D} with
$$\xymatrix@R-.75pc{
M\ar@{ |>->}[d]_{m}\ar@{-->}[rr] & & 
   N\ar@{ |>->}[d]^{n\textstyle{\qquad\mbox{\textit{i.e.}~with}\qquad 
                       m\leq f^{-1}(n).}} \\
X\ar[rr]_-{f} & & Y
}$$

\noindent We shall sometimes refer to this fibration as the ``kernel
fibration''. Every functor $F\colon \Cat{D} \rightarrow \Cat{E}$
in \DCK induces a map of fibrations:
\begin{equation}
\label{KerFibMapEqn}
\vcenter{\xymatrix@R-.75pc{
\KSub(\Cat{D})\ar[rr]\ar[d] & & \KSub(\Cat{E})\ar[d] \\
\Cat{D}\ar[rr]_-{F} & & \Cat{E}
}}
\end{equation}

\noindent because $F$ preserves kernels and pullbacks of kernels---the
latter since pullbacks can be formulated in terms of constructions
that are preserved by $F$, see Lemma~\ref{PullbackLem}. As we shall
see, in some situations, diagram~(\ref{KerFibMapEqn}) is a
pullback---also called a change-of-base situation in this context,
see~\cite{Jacobs99a}. This means that the map $\KSub(X) \rightarrow
\KSub(FX)$ is an isomorphism. 

Let us mention one result about this category $\KSub(\Cat{D})$, which
will be used later.

\begin{lemma}
\label{TotalInvolutionLem}
The category $\KSub(\Cat{D})$ for a dagger category \Cat{D} carries
an involution $\KSub(\Cat{D})\op \rightarrow 
\KSub(\Cat{D})$ given by orthocomplementation:
$$\begin{array}{rclcrcl}
\smash{(M\stackrel{m}{\rightarrowtail} X)}
& \longmapsto &
\smash{(M^{\perp} \stackrel{m\rlap{$^{\scriptscriptstyle\perp}$}}{\rightarrowtail} X)}
& \quad\mbox{and}\quad &
f 
& \longmapsto &
f^{\dag}.
\end{array}$$
\end{lemma}

\begin{proof}
The involution is well-defined because a (necessarily unique) map
$\varphi$ exists if and only if a (necessarily unique) map $\psi$
exists, in commuting squares:
\begin{equation}
\label{InvolutionMapEqn}
\vcenter{\xymatrix@R-1.7pc{
M\ar@{ |>->}[dd]_{m}\ar@{-->}[r]^-{\varphi} & N\ar@{ |>->}[dd]^{n}
& &
\;\;N\rlap{$^{\perp}$}\;\;\ar@{ |>->}[dd]_{n^\perp}\ar@{-->}[r]^-{\psi} 
   & M\rlap{$^\perp$}\ar@{ |>->}[dd]^{m^\perp}
\\
& & \Longleftrightarrow \\
X\ar[r]_-{f} & Y & & Y\ar[r]_-{f^\dag} & X
}}
\end{equation}

\noindent Given $\varphi$, we obtain $\psi$ because $f^{\dag}
\after n^{\perp}$ factors through $\ker(m^{\dag}) = m^{\perp}$ since
$$m^{\dag} \after f^{\dag} \after n^{\perp}
=
\varphi^{\dag} \after n^{\dag} \after n^{\perp}
=
\varphi^{\dag} \after 0
=
0.$$

\noindent The reverse direction follows immediately. 
\end{proof}

\section{Main examples}\label{ExamplesSec}

This section describes our four main examples, namely
\Rel, \PInj, \Hilb and \PHilb, and additionally a general
construction to turn a Boolean algebra into a dagger kernel category.

\subsection{The category \Rel of sets and relations}\label{RelSubsec}

Sets and binary relations $R\subseteq X\times Y$ between them can be
organised in the familiar category \Rel, using relational
composition. Alternatively, such a relation may be described as a
Kleisli map $X\rightarrow \mathcal{P}(Y)$ for the powerset monad
$\mathcal{P}$; in line with this representation we sometimes write
$R(x) = \{y\in Y\,|\, R(x,y)\}$. A third way is to represent such a
morphism in \Rel as (an equivalence class of) a pair of maps $(X
\stackrel{r_1}{\leftarrow} R \stackrel{r_2}{\rightarrow} Y)$ whose
tuple $\langle r_{1}, r_{2}\rangle \colon R\rightarrow X\times Y$ of
legs is injective.

There is a simple dagger on \Rel, given by reversal of
relations: $R^{\dag}(y,x) = R(x,y)$. A map $R\colon X\rightarrow Y$ is
a dagger mono in \Rel if $R^{\dag} \after R = \idmap$, which
amounts to the equivalence:
$$\exin{y}{Y}{R(x,y) \conjun R(x',y)}
\quad\Longleftrightarrow \quad
x=x'$$

\noindent for all $x,x'\in X$. It can be split into two statements:
$$\allin{x}{X}{\exin{y}{Y}{R(x,y)}}
\quad\mbox{and}\quad
\allin{x,x'}{X}{\allin{y}{Y}{R(x,y) \conjun R(x',y) \Rightarrow x=x'}}.$$

\noindent Hence such a dagger mono $R$ is given by a span of the
form
\begin{equation}
\label{RelDagMonoEqn}
\left(\raisebox{1em}{$\xymatrix@C-1.5em@R-1.5em{
& R\ar@{->>}[dl]_{r_{1}}\ar@{ >->}[dr]^{r_{2}} \\
X & & Y
}$}\right)
\end{equation}

\noindent with an surjection as first leg and an injection as second leg. A
dagger epi has the same shape, but with legs exchanged.

The empty set $0$ is a zero object in \Rel, and the resulting
zero map $0\colon X\rightarrow Y$ is the empty relation $\emptyset
\subseteq X\times Y$.

The category \Rel also has kernels. For an arbitrary map $R\colon
X\rightarrow Y$ one takes $\ker(R) = \setin{x}{X}{\neg
  \exin{y}{Y}{R(x,y)}}$ with map $k\colon \ker(R) \rightarrow X$ in
\Rel given by $k(x,x') \Leftrightarrow x=x'$. Clearly, $R \after k =
0$.  And if $S\colon Z\rightarrow X$ satisfies $R\after S = 0$, then
$\neg \exin{x}{X}{R(x,y) \conjun S(z,x)}$, for all $z\in Z$ and $y\in
Y$. This means that $S(z,x)$ implies there is no $y$ with
$R(x,y)$. Hence $S$ factors through the kernel $k$. Kernels are thus
of the following form:
$$\left(\raisebox{1em}{$\xymatrix@C-1.5em@R-1.5em{
& K\ar@{=}[dl]\ar@{ >->}[dr] \\
K & & X
}$}\right) \quad\mbox{with}\quad
K = \{x\in X\,|\, R(x) = \emptyset\}.$$

\auxproof{
We check in detail that \Rel has kernels. First, for $x\in\ker(R)$
and $y\in Y$,
$$\begin{array}{rcl}
(R\after k)(x,y)
& \Longleftrightarrow &
\ex{x'}{k(x,x') \conjun R(x',y)} \\
& \Longleftrightarrow &
\neg\exin{y'}{Y}{R(x,y')} \conjun R(x,y) \\
& \Longleftrightarrow &
\bot \\
& \Longleftrightarrow &
0(x,y).
\end{array}$$

\noindent If $S \after R = 0$, then $S(z,x)$ implies that there
is no $y$ with $R(x,y)$. Hence $S\subseteq Z\times \ker(R)$. 
Moreover, for $z\in Z, x\in X$,
$$\begin{array}{rcl}
(k \after S)(z,x)
& \Longleftrightarrow &
\exin{x'}{\ker(R)}{S(z,x') \conjun k(x',x)} \\
& \Longleftrightarrow &
S(x,z) \conjun \neg\exin{y}{Y}{R(x,y)} \\
& \Longleftrightarrow &
S(x,z).
\end{array}$$
}

\noindent So, kernels are essentially given by subsets: $\KSub(X) =
\mathcal{P}(X)$. Indeed, \Rel is Boolean, in the sense of
Definition~\ref{DagCatKerDef}. A cokernel has the reversed shape.

Finally, a relation $R$ is zero-mono if its kernel is 0, see
Lemma~\ref{ZeroMonoEpiLem}. This means that $R(x)\neq\emptyset$, for
each $x\in X$, so that $R$'s left leg is a surjection.

\begin{proposition}
\label{RelMonosProp}
In \Rel there are proper inclusions:
$$\mbox{kernel}
\;\subsetneq\;
\mbox{dagger mono}
\;\subsetneq\;
\mbox{mono}
\;\subsetneq\;
\mbox{zero-mono}.$$

\noindent Subsets of a set $X$ correspond to kernels in \Rel
with codomain $X$. 
\end{proposition}

There is of course a dual version of this result, for cokernels
and epis.

\begin{proof}
We still need to produce (1)~a zero-mono which is not a mono, and
(2)~a mono which is not a dagger mono. As to~(1), consider $R\subseteq
\{0,1\} \times \{a,b\}$ given by $R = \{(0,a), (1,a)\}$. Its first leg
is surjective, so $R$ is a zero-mono. But it is not a mono: there are
two different relations $\{(*,0)\}, \{(*,1)\} \subseteq \{*\} \times
\{0,1\}$ with $R \after \{(*,0)\} = \{(*,a)\} = R\after \{(*,1)\}$.

As to~(2), consider the relation $R\subseteq \{0,1\}\times\{a,b,c\}$
given by $R = \{(0,a), (0,b),$ $(1,b), (1,c)\}$. Clearly, the first
leg of $R$ is a surjection, and the second one is neither an injection
nor a surjection. We check that $R$ is monic. Suppose $S,T\colon X
\rightarrow \{0,1\}$ satisfy $R\after S = R\after T$. If $S(x,0)$,
then $(R\after S)(x,a) = (R\after T)(x,a)$, so that $T(x,0)$. Similarly,
$S(x,1) \Rightarrow T(x,1)$.  
\end{proof}

We add that the pullback $R^{-1}(n)$ of a kernel $n =
(N=N\rightarrowtail Y)$ along a relation $R\subseteq X\times Y$, as
described in Lemma~\ref{PullbackLem}, is the subset of $X$ given by
the modal formula $\Box_{R}(n)(x) = R^{-1}(n)(x) \Leftrightarrow
(\all{y}{R(x,y) \Rightarrow N(y)})$. As is well-known in modal logic,
$\Box_{R}$ preserves conjunctions, but not
disjunctions. Interestingly, the familiar ``graph'' functor ${\cal
  G}\colon \Sets \rightarrow \Rel$, mapping a set to itself and a
function to its graph relation, yields a map of fibrations
\begin{equation}
\label{SetsRelKerPbEqn}
\vcenter{\xymatrix@R-.75pc{
\Sub(\Sets)\ar[d]\ar[rr] & & \KSub(\Rel)\ar[d] \\
\Sets\ar[rr]_-{\cal G} & & \Rel
}}
\end{equation}

\auxproof{
We check the pullback formulation:
$$\begin{array}{rcl}
R^{-1}(n)
& = &
\ker(\coker(n) \after R) \\
& = &
\ker(\ker(Y\leftarrowtail N=N)^{\dag} \after R) \\
& = &
\ker((\neg N = \neg N \rightarrowtail Y)^{\dag} \after R) \\
& = &
\ker((Y \leftarrowtail \neg N = \neg N) \after R) \\
& = &
\ker(\set{(x,y)}{R(x,y) \conjun \neg N(y)}) \\
& = &
\set{x}{\neg\ex{y}{R(x,y) \conjun \neg N(y)}} \\
& = &
\set{x}{\all{y}{\neg (R(x,y) \conjun \neg N(y))}} \\
& = &
\set{x}{\all{y}{R(x,y) \Rightarrow N(y)}}.
\end{array}$$
}

\noindent which in fact forms a pullback (or a ``change-of-base''
situation, see~\cite{Jacobs99a}). This means that the familiar logic
of sets can be obtained from this kernel logic on relations. In this
diagram we use that inverse image is preserved: for a function
$f\colon X\rightarrow Y$ and predicate $N\subseteq Y$ one has:
$$\begin{array}{rcl}
{\cal G}(f)^{-1}(N)
\hspace*{\arraycolsep} = \hspace*{\arraycolsep}
\Box_{{\cal G}(f)}(N)
& = &
\setin{x}{X}{\all{y}{{\cal G}(f)(x,y) \Rightarrow N(y)}} \\
& = &
\setin{x}{X}{\all{y}{f(x) = y \Rightarrow N(y)}} \\
& = &
\setin{x}{X}{N(f(x))} \\
& = &
f^{-1}(N).
\end{array}$$

\subsection{The category \PInj of sets and partial injections}\label{PInjSubsec}

There is a subcategory \PInj of \Rel also with sets as
objects but with ``partial injections'' as morphisms. These are special
relations $F\subseteq X\times Y$ satisfying $F(x,y) \conjun F(x,y')
\Rightarrow y=y'$ and $F(x,y)\conjun F(x',y)\Rightarrow x=x'$. We shall
therefore often write morphisms $f\colon X\rightarrow Y$ in \PInj 
as spans with the notational convention
$$\left(X \stackrel{f}{\longrightarrow} Y\right) 
= 
\left(\raisebox{1em}{$\xymatrix@C-1.5em@R-1.5em{
& F\ar@{ >->}[dl]_{f_{1}}
           \ar@{ >->}[dr]^{f_{2}} \\
X & & Y
}$}\right),$$

\noindent where spans $\smash{(X \stackrel{f_1}{\leftarrowtail} F
  \stackrel{f_2}{\rightarrowtail} Y)}$ and $\smash{(X
  \stackrel{g_1}{\leftarrowtail} G \stackrel{g_2}{\rightarrowtail}
  Y)}$ are equivalent if there is an isomorphism $\varphi\colon
F\rightarrow G$ with $g_{i} \after \varphi = f_{i}$, for
$i=1,2$---like for relations.

Composition of $\smash{X\stackrel{f}{\rightarrow} Y
  \stackrel{g}{\rightarrow} Z}$ can be described as relational
composition, but also via pullbacks of spans. The identity map
$X\rightarrow X$ is given by the span of identities $X \leftarrowtail
X \rightarrowtail X$. The involution is inherited from \Rel and can be
described as $\smash{\big(X\stackrel{f_1}{\leftarrowtail} F
  \stackrel{f_2}{\rightarrowtail} Y\big)^{\dag}} =
\smash{\big(Y\stackrel{f_2}{\leftarrowtail} F
  \stackrel{f_1}{\rightarrowtail} X\big)}$.

It is not hard to see that $f = \big(X\stackrel{f_1}{\leftarrowtail} F
\stackrel{f_2}{\rightarrowtail} Y\big)$ is a
dagger mono---\textit{i.e.}~satisfies $f^{\dag} \after f =
\idmap$---if and only if its first leg $f_{1}\colon F
\rightarrowtail X$ is an isomorphism. For convenience we therefore
identify a mono/injection $m\colon M \rightarrowtail X$ in
$\Cat{Sets}$ with the corresponding dagger mono
$\big(M\stackrel{\idmap}{\leftarrowtail} M
\stackrel{m}{\rightarrowtail} X\big)$ in \PInj.

\auxproof{The (if)-part is easy, because
it yields a pullback:
$$\xymatrix@C-1em@R-1em{
& & F\ar@{ =}[dl]\ar@{ =}[dr] \\
& F\ar@{ >->}[dl]_{f_{1}}^{\cong}\ar@{ >->}[dr]^{f_{2}} & &
   F\ar@{ >->}[dl]_{f_{2}}\ar@{ >->}[dr]^{f_{1}}_{\cong} \\
X & & Y & & X
}$$

\noindent Conversely, if such an $f$ is a dagger mono then there
is a pullback:
$$\xymatrix@C-1em@R-1em{
& & X\ar@{ >->}[dl]_{e}\ar@{ >->}[dr]^{e} \\
& F\ar@{ >->}[dl]_{f_{1}}\ar@{ >->}[dr]^{f_{2}} & &
   F\ar@{ >->}[dl]_{f_{2}}\ar@{ >->}[dr]^{f_{1}} \\
X & & Y & & X
}$$

\noindent in which $f_{1} \after e = \idmap$. This makes $f_{1}$
an epi, and thus an isomorphism.
}

By duality: $f$ is dagger epi iff $f^{\dag}$ is dagger mono iff the
second leg $f_{2}$ of $f$ is an isomorphism.  Further, $f$ is a
dagger iso iff $f$ is both dagger mono and dagger epi iff both legs
$f_{1}$ and $f_{2}$ of $f$ are isomorphisms.

Like in \Rel, the empty set is a zero object, with corresponding
zero map given by the empty relation, and $0^{\dag} = 0$.

For the description of the kernel of an arbitrary map $\smash{f =
\big(X\stackrel{f_1}{\leftarrowtail} F \stackrel{f_2}{\rightarrowtail} 
Y\big)}$ in \PInj we shall use the \textit{ad hoc} notation
$\neg_{1}F \stackrel{\neg f_1}{\rightarrowtail} X$ for the negation of the
first leg $f_{1}\colon F\rightarrowtail X$, as subobject/subset. It
yields a map: 
$$\ker(f) = \left(\raisebox{1em}{$\xymatrix@C-1.5em@R-1.5em{
& \neg_{1}F\ar@{=}[dl]\ar@{ >->}[dr]^{\neg f_{1}} \\
\neg_{1}F & & X
}$}\right)$$

\noindent It satisfies $f \after \ker(f) = 0$. It is a dagger mono by
construction. Notice that kernels are the same as dagger monos, and
are also the same as zero-monos. They all correspond to subsets, 
so that $\KSub(X) = \mathcal{P}(X)$ and \PInj is Boolean,
like \Rel.

\auxproof{
For an arbitrary map $g\colon Z\rightarrow Y$ with
$f\after g = 0$ we have a pullback:
$$\xymatrix@C-1em@R-1em{
& & 0\ar@{ >->}[dl]\ar@{ >->}[dr] \\
& G\ar@{ >->}[dl]_{g_{1}}\ar@{ >->}[dr]^{g_{2}} & &
   F\ar@{ >->}[dl]_{f_{1}}\ar@{ >->}[dr]^{f_{2}} \\
Z & & X & & Y
}$$

\noindent Hence $g_{2}$ factors through $\neg_{1}F \rightarrowtail X$.
This means that $g$ factors through $\ker(f)$ via a mediating map
$Z\rightarrow \neg_{1}F$.
}

%% The cokernel of $f$ is defined as $\coker(f) = \ker(f^{\dag})^{\dag}$,
%% which is the dagger epi $(Y \stackrel{\neg f_2}{\leftarrowtail}
%% \neg_{2}F \stackrel{\idmap}{\rightarrowtail} \neg_{2}F)$. Notice
%% that for an arbitrary dagger mono $m = (X
%% \stackrel{\idmap}{\leftarrowtail} X \stackrel{m}{\rightarrowtail}
%% Y)$ one has:
%% $$\begin{array}{rcl}
%% \ker(\coker(m))
%% & = &
%% \ker\left(Y \stackrel{\neg m}{\leftarrowtail} \neg M 
%%    \stackrel{\idmap}{\rightarrowtail} \neg M\right) \\
%% & = &
%% \left(\neg\neg M \stackrel{\idmap}{\leftarrowtail} \neg \neg M 
%%    \stackrel{\neg\neg m}{\rightarrowtail} X\right) 
%% \hspace*{\arraycolsep} = \hspace*{\arraycolsep}
%% m.
%% \end{array}$$

The next result summarises what we have seen so far and shows that
\PInj is very different from \Rel (see
Proposition~\ref{RelMonosProp}).

\begin{proposition}
\label{PInjMonosProp}
In \PInj there are proper identities:
$$\mbox{kernel}
\;=\;
\mbox{dagger mono}
\;=\;
\mbox{mono}
\;=\;
\mbox{zero-mono}.$$

\noindent These all correspond to subsets.
\end{proposition}

\subsection{The category \Hilb of Hilbert spaces}\label{HilbSubsec}

Our third example is the category $\Hilb$ of (complex) Hilbert spaces
and continuous linear maps. Recall that a Hilbert space is a vector
space $X$ equipped with an inner product, \textit{i.e.}~a function
$\inprod{-}{-}\colon X \times X \to \field{C}$ that is linear in the
first and anti-linear in the second variable, satisfies $\inprod{x}{x}
\geq 0$ with equality if and only if $x=0$, and $\inprod{x}{y} =
\overline{\inprod{y}{x}}$. Moreover, a Hilbert space must be complete
in the metric induced by the inner product by $d(x,y) =
\sqrt{\inprod{x-y}{x-y}}$.

The Riesz representation theorem provides this category with a dagger.
Explicitly, for $f\colon X \to Y$ a given morphism, $f^\dag
\colon Y \to X$ is the unique morphism satisfying
\[
  \inprod{f(x)}{y}_Y = \inprod{x}{f^\dag(y)}_X
\]
for all $x \in X$ and $y \in Y$.
The zero object is inherited from the category of (complex) vector
spaces: it is the zero-dimensional Hilbert space $\{0\}$, with unique
inner product $\inprod{0}{0}=0$. 

In the category $\Hilb$, dagger mono's are usually called isometries,
because they preserve the metric: $f^\dag \after f = \idmap$ if and only
if
\[
    d(fx,fy) 
  = \inprod{f(x-y)}{f(x-y)}^{\frac{1}{2}}
  = \inprod{x-y}{(f^\dag \after f)(x-y)}^{\frac{1}{2}}
  = d(x,y).
\]
Kernels are inherited from the category of vector spaces. For
$f\colon X \to Y$, we can choose $\ker(f)$ to be (the inclusion of)
$\setin{x}{X}{f(x)=0}$, as this is complete with respect to the
restricted inner product of $X$. Hence kernels correspond to
(inclusions of) closed subspaces. Being inclusions, kernels are
obviously dagger monos. Hence $\Hilb$ is indeed an example of a dagger
kernel category. However, \Hilb is not Boolean. The following
proposition shows that it is indeed different, categorically, from
$\Rel$ and $\PInj$.

\begin{proposition}
\label{HilbMonosProp}
In \Hilb one has:
$$\mbox{kernel}
\;=\;
\mbox{dagger mono}
\;\subsetneq\;
\mbox{mono}
\;=\;
\mbox{zero-mono}.$$
\end{proposition}
\begin{proof}
  For the left equality, notice that both kernels and isometries
  correspond to closed subspaces. It is not hard to show that the
  monos in $\Hilb$ are precisely the injective continuous linear
  functions, establishing the middle proper inclusion. Finally,
  $\Hilb$ has equalisers by $\mathrm{eq}(f,g) = \ker(g-f)$, which
  takes care of the right equality.  
\end{proof}

As is well-known, the $\ell^2$ construction forms a functor $\ell^{2}
\colon \PInj \rightarrow \Hilb$ (but not a functor $\Sets \to \Hilb$), see
\textit{e.g.}~\cite{Barr92,HaghverdiS06}. Since it preserves daggers,
zero object and kernels it is a map in the category \DCK, and therefore
yields a map of kernel fibrations like in~(\ref{KerFibMapEqn}).  It
does not form a pullback (change-of-base) between these fibrations,
since the map $\KSub_{\PInj}(X) = \mathcal{P}(X) \rightarrow
\KSub_{\Hilb}(\ell^2(X))$ is not an isomorphism. 

% \begin{equation}
% \label{pInjHilbKerPbEqn}
% \vcenter{\xymatrix{
% \KSub(\PInj)\ar[d]\ar[rr] & & \KSub(\Hilb)\ar[d] \\
% \PInj\ar[rr]^-{\ell^2} & & \Hilb
% }}
% \end{equation}

\auxproof{

Obviously $\ell^2(\emptyset) = \{0\}$, so that $\ell^2$ preserves zero
objects.  To see that $\ell^2$ preserves kernels, note that
$\ell^2(\ker(f)) \colon \ell^2(X \backslash f_1(Z)) \to \ell^2(X)$ is
determined by $\ell^2(\ker(f))(\varphi)(x) = 0$ if $x \in f_1(Z)$ and
$\ell^2(\ker(f))(\varphi(x)) = \varphi(x)$ if $x \not\in f_1(Z)$.
Also, $\ker(\ell^2(f))$ is (the inclusion of)
  \begin{align*}
      & \setin{\varphi}{\ell^2(X)}{\allin{y}{Y}{\sum_{z \in
           f_2^{-1}(y)} \varphi(f_1(z)) = 0}} \\
    = & \setin{\varphi}{\ell^2(X)}{\allin{z}{Z}{
           \varphi(f_1(z)) = 0}} \\
    = & \setin{\varphi}{\ell^2(X)}{\varphi=0 \mbox{ on } X
           \backslash f_1(Z)},
  \end{align*}
  and is hence isomorphic to $\ell^2(\ker(f))$.

  The following calculation shows that $\ell^2$ also preserves
  $\dag$. For $\varphi \in \ell^2(X)$ and $\psi \in \ell^2(Y)$:
  \begin{align*}
        \inprod{\ell^2(f)(\varphi)}{\psi}_{\ell^2(Y)}
    & = \sum_{y \in Y} \overline{\ell^2(f)(\varphi)(y)} \cdot \psi(y) \\
    & = \sum_{y \in Y} \overline{(\sum_{z \in f_2^{-1}(y)}
                      \varphi(f_1(z)))} \cdot \psi(y) \\
    & = \sum_{y \in Y} \sum_{z \in f_2^{-1}(y)} \overline{\varphi(f_1(z))}
                      \cdot \psi(y) \\
    & = \sum_{z \in Z} \overline{\varphi(f_1(z))} \cdot \psi(f_2(z)) \\
    & = \sum_{x \in X} \sum_{z \in f_1^{-1}(x)} \overline{\varphi(x)}
                      \cdot \psi(f_2(z)) \\
    & = \sum_{x \in X} \overline{\varphi(x)} \cdot 
                      (\sum_{z \in f_1^{-1}(x)} \psi(f_2(z)))\\
    & = \sum_{x \in X} \overline{\varphi(x)} \cdot \ell^2(f^\dag)(\psi)(z) \\
    & = \inprod{\varphi}{\ell^2(f^\dag)(\psi)}_{\ell^2(X)}.
  \end{align*}

}

%% In this paper we concentrate on certain structure in dagger categories
%% and largely ignore structure preserving functors. But it is worthwhile
%% noting that the inclusion and $\ell_2$ functors:
%% $$\xymatrix{
%% \Rel & & \;\PInj\ar@{_{(}->}[ll]\ar[rr]^-{\ell_2} & & \Hilb
%% }$$

%% \noindent preserve all structure from Definitions~\ref{DagCatKerDef}
%% (and also from Definition~\ref{DagCatKerSumDef}).

\subsection{The category \PHilb: Hilbert spaces modulo phase}\label{PHilbSubsec}

  The category $\Cat{PHilb}$ of \emph{projective Hilbert spaces} has
  the same objects as $\Cat{Hilb}$, but its homsets are quotiented by
  the action of the circle group $U(1) =
  \setin{z}{\field{C}}{|z|=1}$. That is, continuous linear
  transformations $f,g\colon X \to Y$ are identified when $x=z\cdot y$
  for some phase $z \in U(1)$.
  
  Equivalently, we could write $PX = X_1 / U(1)$ for an object of
  $\Cat{PHilb}$, where $X \in \Cat{Hilb}$ and $X_1 =
  \setin{x}{X}{\|x\|=1}$. Two vectors $x,y \in X_1$ are therefore
  identified when $x=z \cdot y$ for some $z \in U(1)$. Continuous
  linear transformations $f,g\colon X \to Y$ then descend to the same
  function $PX \to PY$ precisely when they are equivalent under the
  action of $U(1)$. This gives a full functor $P\colon \Hilb \to
  \Cat{PHilb}$.

  The dagger of $\Hilb$ descends to $\Cat{PHilb}$, because
  if $f=z\cdot g$ for some $z \in U(1)$, then
  \[ 
      \inprod{f(x)}{y} 
    = \bar{z} \cdot \inprod{g(x)}{y}
    = \bar{z} \cdot \inprod{x}{g^\dag(y)}
    = \inprod{x}{\bar{z} \cdot g^\dag(y)},
  \]
  whence also $f^\dag = \bar{z} \cdot g^\dag$, making the dagger
  well-defined. 

  Also dagger kernels in $\Cat{Hilb}$ descend to $\Cat{PHilb}$. More
  precisely, the kernel $\ker(f) = \setin{x}{X}{f(x)=0}$ of a morphism
  $f\colon X \to Y$ is well-defined, for if $f=z \cdot f'$ for some $z
  \in U(1)$, then
  \[
      \ker(f) 
    = \setin{x}{X}{z \cdot f'(x)=0}
    = \setin{x}{X}{f'(x)=0}
    = \ker(f').
  \]

\begin{proposition}
\label{PHilbMonosProp}
In \Cat{PHilb} one has:
$$\mbox{kernel}
\;=\;
\mbox{dagger mono}
\;\subsetneq\;
\mbox{mono}
\;=\;
\mbox{zero-mono}.$$
\end{proposition}
\begin{proof}
  It remains to be shown that every zero-mono is a mono. So let
  $m\colon Y \to Z$ be a zero-mono, and $f,g\colon X \to Y$ arbitrary
  morphisms in $\Cat{PHilb}$. More precisely, let $m,f$ and $g$ be
  morphisms in $\Cat{Hilb}$ representing the equivalence classes $[m],
  [f]$ and $[g]$ that are morphisms in $\Cat{PHilb}$. Suppose that $[m
    \after f] = [m \after g]$. Then $m \after f \sim m \after g$, say
  $m \after f = z \cdot (m \after g)$ for $z \in U(1)$. So $m \after
  (f-z \cdot g) = 0$, and $f-z \cdot g=0$ since $m$ is zero-mono. Then
  $f=z \cdot g$ and hence $f \sim g$, \textit{i.e.}~$[f]=[g]$. Thus
  $m$ is mono.  
\end{proof}

The full functor $P\colon\Cat{Hilb} \to \Cat{PHilb}$ preserves
daggers, the zero object and kernels. Hence it is a map in the
category \DCK. In fact it yields a pullback (change-of-base) between the
corresponding kernel fibrations.
\begin{equation}
\label{HilbPHilbKerPbEqn}
\vcenter{\xymatrix@R-.75pc{
\KSub(\Cat{Hilb})\pullback\ar[d]\ar[rr] & & \KSub(\Cat{PHilb})\ar[d] \\
\Cat{Hilb}\ar[rr]_-{P} & & \Cat{PHilb}
}}
\end{equation}

\subsection{From Boolean algebras to dagger kernel categories}\label{BAConstrSubsec}

The previous four examples were concrete categories, to which we add
a generic construction turning an arbitrary Boolean algebra into
a (Boolean) dagger kernel category.

To start, let $B$ with $(1,\conjun)$ be a meet semi-lattice. We can
turn it into a category, for which we use the notation
$\widehat{B}$. The objects of $\widehat{B}$ are elements $x\in B$, and
its morphisms $x\rightarrow y$ are elements $f\in B$ with $f\leq x,y$,
\textit{i.e.}~$f\leq x\conjun y$. There is an identity $x\colon
x\rightarrow x$, and composition of $f\colon x\rightarrow y$ and
$g\colon y\rightarrow z$ is simply $f\conjun g\colon x\rightarrow z$.
This $\widehat{B}$ is a dagger category with $f^{\dag} = f$.  A
map $f\colon x\rightarrow y$ is a dagger mono if $f^{\dag} \after f =
f\conjun f = x$. Hence a dagger mono is of the form $x\colon
x\rightarrow y$ where $x\leq y$.

It is not hard to see that the construction $B\mapsto \widehat{B}$ is
functorial: a morphism $h\colon B\rightarrow C$ of meet semi-lattices
yields a functor $\widehat{h}\colon \widehat{B}\rightarrow \widehat{C}$
by $x\mapsto h(x)$. It clearly preserves $\dag$.

\begin{proposition}
\label{BAConstrProp}
If $B$ is a Boolean algebra, then $\widehat{B}$ is a Boolean dagger
kernel category. This yields a functor $\Cat{BA}\rightarrow
\DCK$.
%\marginpar{Does this need full Boolean structure?} 
\end{proposition}

\begin{proof}
The bottom element $0\in B$ yields a zero object $0\in\widehat{B}$,
and also a zero map $0\colon x\rightarrow y$.  For an arbitrary map
$f\colon x\rightarrow y$ there is a kernel $\ker(f) = \neg f \conjun
x$, which is a dagger mono $\ker(f) \colon \ker(f) \rightarrow x$ in
$\widehat{B}$. Clearly, $f\after \ker(f) = f \conjun \neg f \conjun x
= 0 \conjun x = 0$. If also $g\colon z\rightarrow x$ satisfies
$f\after g = 0$, then $g\leq x,z$ and $f\conjun g = 0$. The latter
yields $g\leq \neg f$ and thus $g\leq \neg f \conjun x = \ker(f)$.
Hence $g$ forms the required mediating map $g\colon z\rightarrow
\ker(f)$ with $\ker(f) \after g = g$.

Notice that each dagger mono $m\colon m\rightarrow x$, where $m\leq
x$, is a kernel, namely of its cokernel $\neg m \conjun x\colon
x\rightarrow (\neg m\conjun x)$.  For two kernels $m\colon
m\rightarrow x$ and $n\colon n\rightarrow x$, where $m,n\leq x$, one
has $m\leq n$ as kernels iff $m\leq n$ in $B$. Thus $\KSub(x) =
\downset x$, which is again a Boolean algebra (with negation
$\neg_{x}m = \neg m \conjun x$). The intersection $m\conjun n$ as
subobjects is the meet $m\conjun n$ in $B$. This allows us to show
that $\widehat{B}$ is Boolean: if $m\conjun n = 0$, them $m^{\dag}
\after n = m \after n = m\conjun n = 0$. 
\end{proof}

\auxproof{
For $m,n\leq x$,
$$\begin{array}{rcl}
\Big(m\stackrel{m}{\rightarrow}x\Big) \leq
   \Big(n\stackrel{n}{\rightarrow}x\Big) 
& \Longleftrightarrow &
\ex{\varphi\leq m, n}{n\conjun \varphi = m} \\
& \Longleftrightarrow &
\ex{\varphi\leq m,n}{n\conjun \varphi = m \mbox{ and } m = n \conjun \varphi \leq \varphi} \\
& \Longleftrightarrow &
m \leq n.
\end{array}$$
}

The straightforward extension of the above construction to
orthomodular lattices does not work: in order to get kernels one needs
to use the and-then connective ($\andthen$, see
Proposition~\ref{SasakiProp}) for composition; but $\andthen$ is
neither associative nor commutative, unless the lattice is
Boolean~\cite{Lehmann08}. However, at the end of~\cite{Jacobs09a} a
dagger kernel category is constructed out of an orthomodular lattice
in a different manner, namely via the (dagger) Karoubi envelope of the
associated Foulis semigroup. For more information about orthomodular
lattices, see~\cite{Kalmbach83}, and for general constructions,
see for instance~\cite{Harding06}.

\section{Factorisation}\label{FactorisationSec}

In this section we assume that \Cat{D} is an arbitrary dagger kernel
category. We will show that each map in \Cat{D} can be factored as 
a zero-epi followed by a kernel, in an essentially unique way. This
factorisation leads to existential quantifiers $\exists$, as is
standard in categorical logic.

The image of a morphism $f\colon X\rightarrow Y$ is defined as
$\ker(\coker(f))$. Since it is defined as a kernel, an image is really 
an equivalence class of morphisms with codomain $X$, up to isomorphism of
the domain. We denote a representing morphism by $i_f$, and its
domain by $\Im(f)$. As with kernels, we can choose $i_f$ to be dagger
mono. Both the morphism $i_f$ and the object $\Im(f)$
are referred to as the image of $f$. Explicitly, it can be obtained in
the following steps. First take the kernel $k$ of $f^\dag$: 
\[\xymatrix{
  \ker(f^{\dag})\ar@{ |>->}[r]^-{k} & Y\ar[r]^-{f^\dag} & X.
}\]
Then define $i_f$ as the kernel of $k^{\dag}$, as in the
following diagram:
\begin{equation}
  \label{ImageEqn}
  \vcenter{\xymatrix{
  \llap{$\Im(f) = \;$}
  \ker(k^{\dag})\ar@{ |>->}[r]^-{i_f} & Y\ar@{-|>}[r]^-{k^{\dag}} & \ker(f^{\dag}). \\
  X\ar[ur]_-{f}\ar@{-->}[u]^{e_f} 
  }}
\end{equation}
The map $e_{f}\colon X\rightarrow \Im(f)$ is obtained from
the universal property of kernels, since $k^{\dag} \after f =
(f^{\dag} \after k)^{\dag} = 0^{\dag} = 0$. Since $i_f$ was chosen
to be dagger mono, this $e_f$ is determined as $e_{f} = \idmap \after
e_{f} = (i_{f})^{\dag} \after i_{f} \after e_{f} = (i_{f})^{\dag} \after f$.

So images are defined as dagger kernels. Conversely, every dagger kernel $m =
\ker(f)$ arises as an image, since $\ker(\coker(m)) = m$ by Lemma~\ref{KerLem}. 

The maps that arise as $e_f$ in~(\ref{ImageEqn}) can be characterised.

\begin{proposition}
\label{FactorisationZeroEpiProp}
The maps in \Cat{D} that arise of the form $e_f$, as in
diagram~(\ref{ImageEqn}), are precisely the zero-epis.
\end{proposition}

\begin{proof}
We first show that $e_f$ is a zero-epi. So, assume a map $h\colon
\ker(k^{\dag}) \rightarrow Z$ satisfying $h \after e_{f} = 0$. Recall
that $e_{f} = (i_{f})^{\dag} \after f$, so that:
$$f^{\dag} \after (i_{f} \after h^{\dag})
=
(h \after (i_{f})^{\dag} \after f)^{\dag}
=
(h \after e_{f})^{\dag} 
=
0^{\dag}
=
0.$$

\noindent This means that $i_{f} \after h^{\dag}$ factors through the
kernel of $f^{\dag}$, say via $a\colon Z\rightarrow \ker(f^{\dag})$
with $k \after a = i_{f} \after h^{\dag}$. Since $k$ is a dagger mono
we now get:
$$a
=
k^{\dag} \after k \after a
=
k^{\dag} \after i_{f} \after h^{\dag}
=
0 \after h^{\dag}
=
0.$$

\noindent But then $i_{f} \after h^{\dag} = k \after a = k \after 0 = 0
= i_{f} \after 0$, so that $h^{\dag} = 0$, because $i_f$ is mono, and
$h = 0$, as required.

Conversely, assume $g\colon X\rightarrow Y$ is a zero-epi, so that
$\coker(g) = 0$ by Lemma~\ref{ZeroMonoEpiLem}. Trivially, $i_{g} =
\ker(\coker(g)) = \ker(X\rightarrow 0) = \idmap[X]$, so that $e_{g} =
g$.  
\end{proof}

The factorisation $f = i_{f} \after e_{f}$ from~(\ref{ImageEqn}) 
describes each map as a zero-epi followed by a kernel.  In fact, these
zero-epis and kernels also satisfy what is usually called the
``diagonal fill-in'' property.

\begin{lemma}
\label{DiagonalFillInLem}
In any commuting square of shape
$$\xymatrix{
\cdot \ar@{->>}|{\circ}[r]\ar[d] & 
   \cdot\ar[d]^{\qquad\mbox{there is a (unique) diagonal}} 
& \hspace*{10em} & 
\cdot \ar@{->>}|{\circ}[r]\ar[d] & 
   \cdot\ar[d]\ar@{-->}[dl] \\
\cdot\ar@{ |>->}[r] & \cdot
& &
\cdot\ar@{ |>->}[r] & \cdot
}$$

\noindent making both triangles commute.

As a result, the factorisation~(\ref{ImageEqn}) is unique
up to isomorphism. Indeed, kernels and zero-epis form a factorisation
system (see~\cite{BarrW85}).
\end{lemma}

\begin{proof}
Assume the zero-epi $e\colon E\rightarrow Y$ and kernel $m = \ker(h)
\colon M\rightarrowtail X$ satisfy $m\after f = g \after e$, as
below,
$$\xymatrix@C-.5pc@R-.5pc{
E\ar@{->>}|{\circ}[r]^-{e}\ar[d]_{f} & Y\ar[d]^{g} \\
M\ar@{ |>->}[r]_-{m} & X\ar[r]_-{h} & Z
}$$

\noindent Then: $h\after g \after e = h \after m \after f = 0 \after f
= 0$ and $h\after g = 0$ because $e$ is zero-epi. This yields the
required diagonal $d\colon Y\rightarrow M$ with $m \after d = g$
because $m$ is the kernel of $h$. Using that $m$ is monic we get $d
\after e = f$.  
\end{proof}

Factorisation standardly gives a left adjoint to inverse image
(pullback), corresponding to existential quantification in logic. In
this self-dual situation there are alternative descriptions.

Notice that this general prescription of quantifiers by categorical
logic, when applied to our quantum setting, is of a different 
nature from earlier attempts at quantifiers for quantum
logic~\cite{Janowitz63,Roman06}, as it 
concerns multiple orthomodular lattices instead of a single one.

\begin{proposition}
\label{ExistentialProp}
For $f\colon X\rightarrow Y$, the pullback functor $f^{-1}\colon
\KSub(Y)\rightarrow \KSub(X)$ from Lemma~\ref{PullbackLem} has
a left adjoint $\exists_f$ given as image:
$$\xymatrix{
\Big(M\ar@{ |>->}[r]^-{m} & X\Big) \;\longmapsto \; 
   \Big(\Im(f\after m) \ar@{ |>->}[rr]^-{\exists_{f}(m)=i_{f\after m}}
   & & Y\Big)
}$$

\noindent Alternatively, $\exists_{f}(m) =
\Big((f^{\dag})^{-1}(m^{\perp})\Big)^{\perp}$.
\end{proposition}
\begin{proof}
  The heart of the matter is that in the following diagram, the map
  $\varphi$ (uniquely) exists if and only if the map $\psi$ (uniquely)
  exists: 
  \[\xymatrix@C+3ex{
    M \ar@{ |>->}@(d,l)_-{m}[dr] \ar@{-->}^-{\varphi}[r]
    \ar@{-->}@(ur,ul)^-{\psi}[rr]  
    & \cdot \ar@{ |>->}_-{f^{-1}(n)}[d] \pullback \ar[r] 
    & N \ar@{ |>->}^-{n}[d] \\
    & X \ar_-{f}[r] & Y\rlap{.}
  }\]
  Thus one easily reads off:
  \begin{align*}
      m \leq f^{-1}(n)
    & \Longleftrightarrow \text{ there is } \varphi \text{ such that }
    m=f^{-1}(n) \after \varphi \\
    & \Longleftrightarrow \text{ there is } \psi \text{ such that } f
    \after m = n \after \psi \\
    & \Longleftrightarrow \exists_f(m) \leq n.
  \end{align*}
  For the alternative description:
$$\begin{array}[b]{rcccl}
\Big((f^{\dag})^{-1}(m^{\perp})\Big)^{\perp} \leq n
& \Longleftrightarrow &
n^{\perp} \leq (f^{\dag})^{-1}(m^{\perp}) 
& \smash{\stackrel{(\ref{InvolutionMapEqn})}{\Longleftrightarrow}} &
m \leq f^{-1}(n).
\end{array}$$
\end{proof}

This adjunction $\exists_{f} \dashv f^{-1}$ makes the kernel fibration
$\Big(\fibno{\KSub(\Cat{D})}{\Cat{D}}\Big)$ an opfibration, and thus a
bifibration, see~\cite{Jacobs99a}.
Recall the \emph{Beck-Chevalley condition}: if the left square
below is a pullback in $\Cat{D}$, then the right one must commute.
\begin{equation}
\label{eq:beckchevalley}\tag{BC}
  \raise5ex\hbox{\xymatrix{
    P \pullback \ar^-{q}[r] \ar_-{p}[d] & Y \ar^-{g}[d] \\
    X \ar_-{f}[r] & Z    
  }} \qquad \Longrightarrow \qquad
  \raise5ex\hbox{\xymatrix{
    \KSub(P) \ar_-{\exists_p}[d] & \KSub(Y) \ar^-{\exists_g}[d]
    \ar_-{q^{-1}}[l] \\
    \KSub(X) & \KSub(Z) \ar^-{f^{-1}}[l] 
  }}
\end{equation}
This condition ensures that $\exists$ commutes with substitution. If
one restricts attention to the pullbacks of the form given in
Lemma~\ref{PullbackLem}, then Beck-Chevalley holds. In the notation of
Lemma~\ref{PullbackLem}, for kernels $k\colon K\rightarrowtail Y$ and
$g\colon Y \rightarrowtail Z$:
$$\begin{array}{rcll}
f^{-1}(\exists_g(k))
& = &
f^{-1}(g \after k)
   & \mbox{because both $g,k$ are kernels} \\
& = &
p \after q^{-1}(k) 
   & \mbox{by composition of pullbacks} \\
& = &
\exists_p(q^{-1}(k)).
\end{array}$$

\noindent In \Hilb all pullbacks exist and Beck-Chevalley holds for
all of them by~\cite[II, Proposition~1.7.6]{Borceux94} using \Hilb's
biproducts and equalisers.

%% \begin{lemma}
%% \label{DSubPullbackLem}
%% Kernels are stable under pullback.
%% \end{lemma}

%% \begin{proof}
%% Assume a pullback as below on the left, in which we may already assume
%% that $m$ is monic (because monos are stable under pullback).  On the
%% right we factor $m$ as $m = i_{m}\after e_{m}$, like
%% in~(\ref{ImageEqn}).
%% $$\xymatrix{
%% M\pullback \ar[rr]^-{f'}\ar@{ >->}[d]_{m} & & N\ar@{ |>->}[d]^{n = \ker(h)}
%% & &
%% M\pullback \ar[rr]^-{f'}\ar@{->>}|{\circ}[d]_{e_m} & & 
%%    N\ar@{ |>->}[dd]^{n = \ker(h)} \\
%% X\ar[rr]_-{f} & & Y\ar@{-|>}[d]^{h} 
%% & &
%% \Im(m)\ar@{ |>->}[d]_{i_m}\ar@{-->}[urr]_-{g} \\
%% & & & &
%% X\ar[rr]_-{f} & & Y 
%% }$$

%% \noindent Then: $h \after f \after i_{m} \after e_{m} = h \after n
%% \after f' = 0 \after f' = 0$, so that $h\after f \after i_{m} = 0$
%% because $e_m$ is zero-epi. This means that there is a map $g\colon
%% \Im(m) \rightarrow N$ with $n \after g = f \after i_{m}$, as
%% indicated. The pullback property then yields a map $\varphi \colon
%% \Im(m) \rightarrow M$ with $f' \after \varphi = g$ and $m \after
%% \varphi = i_{m}$. Then:
%% $$\begin{array}{rcl}
%% e_{m} \after \varphi = \idmap
%% & \mbox{because} &
%% i_{m} \after e_{m} \after \varphi = m \after \varphi = i_{m} \\
%% \varphi \after e_{m} = \idmap
%% & \mbox{because} &
%% m \after \varphi \after e_{m} = i_{m} \after e_{m} = m.
%% \end{array}$$

%% \noindent Hence $e_{m}$ is an iso and $m$ is isomorphic
%% to a kernel, and hence a kernel itself. 
%% \end{proof}

The final result in this section brings more clarity; it underlies the
relations between the various maps in the propositions in the previous
section.

\begin{lemma}
\label{ZeroEpiEpiLem}
If zero-epis are (ordinary) epis, then dagger monos are kernels.
\end{lemma}

Recall that Lemma~\ref{EqualiserZeroEpiLem} tells that zero-epis
are epis in the presence of equalisers.

\begin{proof}
Suppose $m\colon M\rightarrowtail X$ is a dagger mono, with
factorisation $m = i \after e$ as in~(\ref{ImageEqn}), where $i$ is a
kernel and a dagger mono, and $e$ is a zero-epi and hence an epi by
assumption. We are done if we can show that $e$ is an isomorphism.
Since $m = i \after e$ and $i$ is dagger monic we get $i^{\dag} \after
m = i^{\dag} \after i \after e = e$. Hence $e^{\dag} \after e =
(i^{\dag} \after m)^{\dag} \after e = m^{\dag} \after i \after e =
m^{\dag} \after m = \idmap$ because $m$ is dagger mono. But then also
$e \after e^{\dag} = \idmap$ because $e$ is epi and $e\after e^{\dag}
\after e = e$.  
\end{proof}

\begin{example}
\label{ImageDomainEx}
In the category \Rel the image of a morphism $(X
\stackrel{r_1}{\leftarrow} R \stackrel{r_2}{\rightarrow} Y)$ is the
relation $\smash{i_{R} = (Y' \stackrel{=}{\leftarrow} Y'
  \rightarrowtail Y)}$ where $Y' = \setin{y}{Y}{\ex{x}{R(x,y)}}$ is
the image of the second leg $r_2$ in \Cat{Sets}. The associated
zero-epi is $e_{R} = 
(X \stackrel{r_1}{\leftarrow} R \stackrel{r_2}{\twoheadrightarrow}
Y')$. Existential quantification $\exists_{R}(M)$ from
Proposition~\ref{ExistentialProp} corresponds to the modal diamond
operator (for the reversed relation $R^{\dag}$):
$$\exists_{R}(M)
=
\setin{y}{Y}{\exin{x}{M}{R(x,y)}}
=
\diamondsuit_{R^{\dag}}(M)
=
\neg\Box_{R^{\dag}}(\neg M).$$

\auxproof{
We check the adjunction:
$$\begin{array}{rcl}
\exists_{R}(M) \subseteq N
& \Longleftrightarrow &
\all{y}{(\ex{x}{R(x,y)\conjun M(x)}) \Rightarrow N(y)} \\
& \Longleftrightarrow &
\all{x,y}{R(x,y)\conjun M(x) \Rightarrow N(y)} \\
& \Longleftrightarrow &
\all{x,y}{M(x) \Rightarrow (R(x,y) \Rightarrow N(y))} \\
& \Longleftrightarrow &
\all{x}{M(x) \Rightarrow \all{y}{R(x,y) \Rightarrow N(y)}} \\
& \Longleftrightarrow &
M \subseteq R^{-1}(N).
\end{array}$$
}

\noindent It is worth mentioning that the ``graph'' map of
fibrations~(\ref{SetsRelKerPbEqn}) between sets and relations is
also a map of opfibrations: for a function $f\colon X\rightarrow Y$
and a predicate $M\subseteq X$ one has:
$$\begin{array}{rcl}
\exists_{{\cal G}(f)}(M)
& = &
\set{y}{\ex{x}{{\cal G}(f)(x,y) \conjun M(x)}} \\
& = &
\set{y}{\ex{x}{f(x) = y \conjun M(x)}} \\
& = &
\set{f(x)}{M(x)} \\
& = &
\exists_{f}(M),
\end{array}$$

\noindent where $\exists_f$ in the last line is the left adjoint
to pullback $f^{-1}$ in the category \Sets.

In \PInj the image of a map $f = (X \stackrel{f_1}{\leftarrowtail} F
\stackrel{f_2}{\rightarrowtail} Y)$ is given as $i_{f} = (F
\stackrel{\idmap}{\leftarrowtail} F \stackrel{f_2}{\rightarrowtail}
Y)$. The associated map $e_f$ is $(X \stackrel{f_1}{\leftarrowtail} F
\stackrel{\idmap}{\rightarrowtail} F)$, so that indeed $f = i_{f}
\after e_{f}$. Notice that this $e_f$ is a dagger epi in \PInj.

In \Hilb, the image of a map $f\colon X \to Y$ is (the inclusion of)
the closure of the set-theoretic image
$\setin{y}{Y}{\exin{x}{X}{y=f(x)}}$. This descends to $\Cat{PHilb}$:
the image of a morphism is the equivalence class represented by the
inclusion of the closure of the set-theoretic image of a
representative.

The functor $\ell^2 \colon \PInj \to \Hilb$ is a map of opfibrations:
for a partial injection $f = (X \stackrel{f_1}{\leftarrowtail} F
\stackrel{f_2}{\rightarrowtail} Y)$ and a kernel $m\colon M
\rightarrowtail X$ in $\PInj$ one has:
\begin{align*}
      \exists_{\ell^2(f)}(\ell^2(m)) 
  & = \mathrm{Im}_{\Hilb}(\ell^2(f \after m)) \\
  & = \mathrm{Im}_{\Hilb}(\ell^2(M) \times Y \ni (\varphi,y) \mapsto \sum_{x
        \in (f \after m)^{-1}(y)} \varphi(x))) \\
  & \cong \overline{\setin{\varphi}{\ell^2(X)}{\mathrm{supp}(\varphi) \subseteq
         F \cap M}} \\
  & = \setin{\varphi}{\ell^2(X)}{\mathrm{supp}(\varphi) \subseteq
         F \cap M} \\
%  & = \ell^2(\mathrm{dom}(f) \cap M) \\
  & \cong \ell^2(f_2 \after f_1^{-1}(m)) \\
  & = \ell^2(\exists_f(m)).
\end{align*}

Also the full functor $P\colon \Hilb \to \PHilb$ is a map of
opfibrations: for $f\colon X \to Y$ and a kernel $m\colon M
\rightarrowtail X$ in $\Hilb$ one has:
\begin{align*}
      \exists_{Pf}(Pm)
  & = \mathrm{Im}_{\PHilb}(P(f \after m)) \\
  & = \overline{\set{f(x)}{x \in M}} \\
  & = P(\overline{\set{f(x)}{x \in M}}) \\
  & = P(\mathrm{Im}_{\Hilb}(f \after m)) \\
  & = P(\exists_f(m)).
\end{align*}

In the category $\widehat{B}$ obtained from a Boolean algebra the
factorisation of $f\colon x\rightarrow y$ is the composite $x
\smash{\stackrel{f}{\longrightarrow}} f
\smash{\stackrel{f}{\longrightarrow}} y$. In particular, for $m\leq x$,
considered as kernel $m\colon m\rightarrow x$ one has $\exists_{f}(m) = 
(m\conjun f\colon (m\conjun f)\rightarrow x)$.
\end{example}

\begin{example}
\label{ManesEx}
In~\cite{Manes89} the domain $\Dom(f)$ of a map $f\colon X\rightarrow
Y$ is the complement of its kernel, so $\Dom(f) = \ker(f)^{\perp}$,
and hence a kernel itself. It can be described as an image, namely of
$f^{\dag}$, since:
$$\Dom(f)
=
\ker(f)^{\perp}
=
\ker(\ker(f)^{\dag})
=
\ker(\coker(f^{\dag}))
=
i_{f^{\dag}}.$$

\noindent It is shown in~\cite{Manes89} that the composition $f
\after \Dom(f)$ is zero-monic---or ``total'', as it is called there.
This also holds in the present setting, since:
$$f \after \Dom(f)
=
f^{\dag\dag} \after i_{f^{\dag}}
=
(i_{f^{\dag}} \after e_{f^{\dag}})^{\dag} \after i_{f^{\dag}}
=
(e_{f^{\dag}})^{\dag} \after (i_{f^{\dag}})^{\dag} \after i_{f^{\dag}}
=
(e_{f^{\dag}})^{\dag}.$$

\noindent This $e_{f^{\dag}}$ is zero-epic, by
Proposition~\ref{FactorisationZeroEpiProp}, so that
$(e_{f^{\dag}})^{\dag}$ is indeed zero-monic. In case $f\colon
X\rightarrow X$ is a self-adjoint map, meaning $f^{\dag} = f$, then
the image of $f$ is the same as the domain, and thus as the 
complement of the kernel.
\end{example}

There is one further property that is worth making explicit, if only
in examples. In the kernel fibration over \Rel one finds the following
correspondences.
$$\KSub(X)
\cong
{\cal P}(X)
\cong
\Sets(X, 2)
\cong
\Sets(X, {\cal P}(1))
\cong
\Rel(X, 1).$$
This suggests that one has ``kernel classifiers'',
comparable to ``subobject classifiers'' in a topos---or more
abstractly, ``generic objects'', see~\cite{Jacobs99a}. But the
naturality that one has in toposes via pullback functors $f^{-1}$
exists here via their left adjoints $\exists_f$. That is, we really
have found a natural correspondence $\KSub(X) \cong \Rel(1,X)$ instead of
$\KSub(X) \cong \Rel(X,1)$. Indeed, there are
natural ``characteristic'' isomorphisms:
$$\xymatrix@R-2pc{
\llap{$\KSub(X)=\;$}{\cal P}(X)\ar[rr]^-{\charac}_-{\cong} & & \Rel(1,X) \\
(M\subseteq X)\ar@{|->}[rr] & & \set{(*,x)}{x\in M}.
}$$

\noindent Then, for $S\colon X\rightarrow Y$ in \Rel,
$$\begin{array}{rcl}
S \after \charac(M)
& = &
\set{(*,y)}{\ex{x}{\charac(M)(*,x) \conjun S(x,y)}} \\
& = &
\set{(*,y)}{\ex{x}{M(x) \conjun S(x,y)}} \\
& = &
\set{(*,y)}{\exists_{S}(M)(y)} \\
& = &
\charac(\exists_{S}(M)).
\end{array}$$
Hence one could say that \Rel has a kernel ``opclassifier''.
This naturality explains our choice of $\Rel(1,X)$ over $\Rel(X,1)$:
the latter formulation more closely resembles the subobject
classifiers of a topos, but using the former, naturality can be
formulated without using the dagger. Hence in principle one could even
consider ``opclassifiers'' in categories without a dagger.

The same thing happens in the dagger categories $\widehat{B}$ from
Subsection~\ref{BAConstrSubsec}. There one has, for $x\in B$,
$$\xymatrix@R-2pc{
\llap{$\KSub(x)=\;$}\downset x\ar[rr]^-{\charac}_-{\cong} & & 
   \widehat{B}(1,x) \\
(m\leq x)\ar@{|->}[rr] & & (m\colon 1\rightarrow x)
}$$

\noindent As before, $f \after \charac(m) = f\conjun m =
\exists_{f}(m) = \charac(\exists_{f}(m))$. 

The category \Cat{OMLatGal} of orthomodular lattices and Galois
connections between them from~\cite{Jacobs09a} also has such an
opclassifier.
There is no obvious kernel opclassifier for the category
$\Cat{Hilb}$. The category $\Cat{PInj}$ is easily seen not to have a
kernel opclassifier.

\section{Images and coimages}\label{ImageCoimageSec}

We continue to work in an arbitrary dagger kernel category \Cat{D}.
In the previous section we have seen how each map $f\colon
X\rightarrow Y$ in \Cat{D} can be factored as $f = i_{f} \after e_{f}$
where the image $i_{f} = \ker(\coker(f)) \colon\Im(f)\rightarrowtail
Y$ is a kernel and $e_f$ is a zero-epi. We can apply this same
factorisation to the dual $f^{\dag}$. The dual of its image,
$(i_{f^\dag})^{\dag} = \coker(\ker(f)) \colon X \twoheadrightarrow
\Im(f^{\dag})$, is commonly called the coimage of $f$.  It is a
cokernel and dagger epi by construction. Thus we have:
\[\xymatrix@R-3ex@C+1ex{
X\ar@{->}[rr]^-{f}\ar@{->>}|{\circ}[dr]_{e_{f}} & & Y
& & 
Y\ar@{->}[rr]^-{f^\dag}\ar@{->>}|{\circ}[dr]_{e_{f^\dag}} & & X \\
& \Im(f)\ar@{ |>->}[ur]_{i_{f}} &
& &
& \Im(f^{\dag})\ar@{ |>->}[ur]_{i_{f^\dag}} &
}\]

\noindent By combining these factorisations we get two
mediating maps $m$ by diagonal fill-in (see
Lemma~\ref{DiagonalFillInLem}), as in: 
\[\xymatrix@R-1ex@C+2ex{
X\ar@{->}[rr]^-{f}\ar@{->>}|{\circ}[dr]^{e_{f}}
   \ar@{-|>}[ddr]_{(i_{f^\dag})^{\dag}} & & Y
& 
Y\ar@{->}[rr]^-{f^\dag}\ar@{->>}|{\circ}[dr]^{e_{f^\dag}}
   \ar@{-|>}[ddr]_{(i_{f})^{\dag}} & & X \\
& \Im(f)\ar@{ |>->}[ur]^{i_{f}} &
& 
& \Im(f^{\dag})\ar@{ |>->}[ur]^{i_{f^\dag}} & \\
& \Im(f^{\dag})\ar@{ >->}|{\circ}[uur]_{(e_{f^\dag})^{\dag}}\ar@{..>}[u]|{m_f} &
& 
& \Im(f)\ar@{ >->}|{\circ}[uur]_{(e_{f})^{\dag}}\ar@{..>}[u]|{m_{f^\dag}} & 
}\]

\noindent We claim that $(m_{f})^{\dag} = m_{f^{\dag}}$. This 
follows easily from the fact that $(i_{f^\dag})^{\dag}$ is epi:
\[ (m_{f^{\dag}})^{\dag} \after (i_{f^\dag})^{\dag}
=
(i_{f^\dag} \after m_{f^{\dag}})^{\dag} 
=
(e_{f})^{\dag\dag} 
=
e_{f}
=
m_{f} \after (i_{f^\dag})^{\dag}. \]

\noindent Moreover, $m_f$ is both a zero-epi and a zero-mono.

\auxproof{ 
If $g \after m_f = 0$, then also $g \after i_f^\dag \after
(e_{f^\dag})^\dag = g \after i_f^\dag \after i_f \after m_f = g
\after m_f = 0$, so $g \after i_f^\dag = 0 = 0 \after i_f^\dag$,
whence $g=0$.

Thus $m_{f^{\dag}}$ is also zero-epi. Hence $m_{f} = (m_{f^{\dag}})^{\dag}$
is zero-mono. Explicitly,

It is also zero-mono: if $m_f \after g = 0$, then $g^\dag \after
e_{f^\dag} \after i_f^\dag = g^\dag \after m_{f^\dag} \after i_f
\after i_f^\dag = g^\dag \after m_f^\dag = 0$, so $g^\dag \after
e_{f^\dag} = 0$ whence $g^\dag=0$ and $g=0$.
}

As a result we can factorise each map $f\colon X\rightarrow Y$ 
in \Cat{D} as:
\begin{equation}
\label{ImageCoimageEqn}
\vcenter{\xymatrix{
X\ar@{-|>}[rr]^-{(i_{f^\dag})^{\dag}}_-{\textrm{coimage}}
& & \Im(f^{\dag})\ar@{ >->>}[r]|{\circ}^-{m_f}_{\begin{array}{c}
   \labelstyle\textrm{zero-epi} \\[-.7pc] 
   \labelstyle\textrm{zero-mono} \end{array}}
& \Im(f)\ar@{ |>->}[rr]^-{i_f}_-{\textrm{image}} & & Y.
}}
\end{equation}

\noindent This coimage may also be reversed, so that a map
in \Cat{D} can also be understood as a pair of kernels
with a zero-mono/epi between them, as in:
$$\xymatrix{
X
& & \Im(f^{\dag})\ar@{ >->>}[r]|{\circ}\ar@{ |>->}[ll]_-{i_{f^\dag}}
& \Im(f)\ar@{ |>->}[rr]^-{i_f} & & Y
}$$

\noindent The two outer kernel maps perform some ``bookkeeping'' to
adjust the types; the real action takes place in the middle, see the
examples below. The category \PInj consists, in a sense, of only
these bookkeeping maps, without any action. This will be described
more systematically in Definition~\ref{KcKDef}.

\begin{example}
\label{ImageCoimageEx}
We briefly describe the factorisation~(\ref{ImageCoimageEqn})
in \Rel, \PInj and \Hilb, using diagrammatic order
for convenience (with notation $f;g = g \after f$).

For a map $(X \stackrel{r_1}{\leftarrow} R \stackrel{r_2}{\rightarrow}
Y)$ in \Rel we take the images $X' \rightarrowtail X$ of $r_1$
and $Y'\rightarrowtail Y$ of $r_2$ in:
$$\left(\raisebox{1em}{$\xymatrix@C-2em@R-1.5em{
& R\ar[dl]_{r_1}\ar[dr]^{r_2} \\
X & & Y
}$}\right)
=
\left(\raisebox{1em}{$\xymatrix@C-2em@R-1.5em{
& X'\ar@{ >->}[dl]\ar@{=}[dr] \\
X & & X'
}$}\right)
;
\left(\raisebox{1em}{$\xymatrix@C-2em@R-1.5em{
& R\ar@{->>}[dl]_{r_1}\ar@{->>}[dr]^{r_2} \\
X' & & Y'
}$}\right)
;
\left(\raisebox{1em}{$\xymatrix@C-2em@R-1.5em{
& Y'\ar@{=}[dl]\ar@{ >->}[dr] \\
Y' & & Y
}$}\right)$$

In \PInj the situation is simpler, because the middle part
$m$ in~(\ref{ImageCoimageEqn}) is the identity, in:
$$\left(\raisebox{1em}{$\xymatrix@C-2em@R-1.5em{
& F\ar@{ >->}[dl]_{f_1}\ar@{ >->}[dr]^{f_2} \\
X & & Y
}$}\right)
=
\left(\raisebox{1em}{$\xymatrix@C-2em@R-1.5em{
& F\ar@{ >->}[dl]_{f_1}\ar@{=}[dr] \\
X & & F
}$}\right)
;
\left(\raisebox{1em}{$\xymatrix@C-2em@R-1.5em{
& F\ar@{=}[dl]\ar@{ >->}[dr]^{f_2} \\
F & & Y
}$}\right).$$

In \Hilb, a morphism $f\colon X \to Y$ factors as $f=i \after m \after
e$. The third part $i\colon I \to Y$ is given by $i(y)=y$, where $I$
is the closure $\overline{\{f(x) \;\colon x \in X\}}$. The first part
$e\colon X \to E$ is given by orthogonal projection on the closure $E
= \overline{\{f^\dag(y) \;\colon y \in Y\}}$; explicitly, $e(x)$ is
the unique $x'$ such that $x=x'+x''$ with $x' \in E$ and
$\inprod{x''}{z}=0$ for all $z \in E$.  Using the fact that the
adjoint $e^\dag \colon E \to X$ is given by $e^\dag(x)=x$, we deduce
that the middle part $m\colon E \to I$ is determined by $m(x) = (i
\after m)(x) = (f \after e^\dag)(x) = f(x)$.  Explicitly,
$$\left(X\stackrel{f}{\longrightarrow} Y\right)
=
\left(X\stackrel{e}{\longrightarrow} E\right)
;
\left(E\stackrel{m}{\longrightarrow} I\right)
;
\left(I\stackrel{i}{\longrightarrow} Y\right).$$
\end{example}

\section{Categorical logic}\label{LogicSec}

This section further investigates the logic of dagger kernel
categories. We shall first see how the so-called Sasaki
hook~\cite{Kalmbach83} arises naturally in this setting, and then
investigate Booleanness.

For a kernel $m\colon M\rightarrowtail X$ we shall write $\effect{m} =
m \after m^{\dag} \colon X\rightarrow X$ for the ``effect'' of $m$,
see~\cite{DvurecenskijP00}. This $\effect{m}$ is easily seen to be a
self-adjoint idempotent: one has $\effect{m}^{\dag} = \effect{m}$ and
$\effect{m} \after \effect{m} = \effect{m}$.  The endomap
$\effect{m}\colon X\rightarrow X$ associated with a kernel/predicate
$m$ on $X$ maps everything in $X$ that is in $m$ to itself, and what
is perpendicular to $m$ to $0$, as expressed by the equations
$\effect{m} \after m = m$ and $\effect{m} \after m^{\perp} = 0$. Of
interest is the following result. It makes the dynamical aspects of
quantum logic described in~\cite{CoeckeS04} explicit.

\begin{proposition}
\label{SasakiProp}
For kernels $m\colon M\rightarrowtail X$, $n\colon N\rightarrowtail X$
the pullback $\effect{m}^{-1}(n)$ is the Sasaki hook, written here as
$\sasaki$:
$$m \sasaki n 
\smash{\;\stackrel{\textrm{def}}{=}\;} 
\effect{m}^{-1}(n)
\;=\;
m^{\perp} \disjun (m\conjun n).$$

\noindent The associated left adjoint $\exists_{\effect{m}} \dashv
\effect{m}^{-1}$ yields the ``and then'' operator:
$$k \andthen m
\smash{\;\stackrel{\textrm{def}}{=}\;} 
\exists_{\effect{m}}(k)
\;=\;
m \conjun (m^{\perp} \disjun k),$$

\noindent so that the ``Sasaki adjunction'' (see~\cite{Finch70}) holds
by construction:
$$\begin{array}{rcl}
k \andthen m \leq n
& \Longleftrightarrow &
k \leq m\sasaki n.
\end{array}$$
\end{proposition}

Quantum logic based on this ``and-then'' $\andthen$ connective is
developed in~\cite{Lehmann08}, see also~\cite{RomanR91,RomanZ96}.
This $\andthen$ connective is in general non-commutative and
non-associative\footnote{The ``and-then'' connective $\andthen$ should
not be confused with the multiplication of a
quantale~\cite{Rosenthal90}, since the latter is always
associative.}. Some basic properties are: $m\andthen m = m$, 
$1\andthen m = m\andthen 1 = m$, $0\andthen m = m\andthen 0 = 0$, and
both $k\andthen m\leq n$, $k^{\perp} \andthen m \leq n$ imply $m\leq
n$ (which easily follows from the Sasaki adjunction).

\begin{proof}
Consider the following pullbacks.
$$\xymatrix{
P\ar[d]_{p}\ar[r]^-{q}\pullback & N\ar@{ |>->}[d]^{n}
& & 
Q\ar[d]_{r}\ar[r]^-{s}\pullback & 
   P\rlap{$^{\perp}$}\ar@{ |>->}[d]^{(m\conjun n)^{\perp} =\,\ker(p^{\dag} \after m^{\dag})} \\
M\ar@{ |>->}[r]_-{m} & X
& &
M\ar@{ |>->}[r]_-{m} & X
}$$

\noindent Then:
$$\begin{array}[b]{rcl}
m^{\perp} \disjun (m\conjun n)
& = &
\big(m \conjun (m\conjun n)^{\perp}\big)^{\perp} \\
& = &
\ker\big((m \conjun (m\conjun n)^{\perp})^{\dag}\big) \\
& = &
\ker\big(r^{\dag} \after m^{\dag}\big) \\
& = &
\ker\big(\ker(\coker((m\conjun n)^{\perp}) \after m)^{\dag} \after m^{\dag}\big) \\
& & \qquad\mbox{by definition of $r$ as pullback, 
    see Lemma~\ref{PullbackLem}} \\
& = &
\ker\big(\ker(\coker(\ker(p^{\dag} \after m^{\dag})) \after m)^{\dag} \after 
  m^{\dag}\big) \\
& = &
\ker\big(\ker(p^{\dag} \after m^{\dag} \after m)^{\dag} \after 
  m^{\dag}\big) \\
& & \qquad\mbox{because $p^{\dag} \after m^{\dag}$ is a cokernel,
   see Lemma~\ref{KernelCompLem}} \\
& = &
\ker(\coker(p) \after m^{\dag}) \\
& = &
\big(m^{\dag}\big)^{-1}(p) \\
& = &
\big(m^{\dag}\big)^{-1}(m^{-1}(n)) \\
& = &
\effect{m}^{-1}(n).
\end{array}$$
\end{proof}

\auxproof{
We check the equivalence $k \andthen m \leq n
\Longleftrightarrow k \leq m\sasaki n$ explicitly.

$(\Rightarrow)$ If $k \andthen m \leq n$ then:
$$\begin{array}{rcl}
m\sasaki n
& = &
m^{\perp} \disjun (m\conjun n) \\
& \geq &
m^{\perp} \disjun (m\conjun (k \andthen m)) \\
& = &
m^{\perp} \disjun (m\conjun (m \conjun (m^{\perp} \disjun k))) \\
& = &
m^{\perp} \disjun (m \conjun (m^{\perp} \disjun k)) \\
& = &
m^{\perp} \disjun k \qquad 
   \mbox{by orthomodularity, since }m^{\perp} \disjun k \geq m^{\perp} \\
& \geq &
k.
\end{array}$$

$(\Leftarrow)$ If $k \leq m\sasaki n$, then:
$$\begin{array}{rcl}
k \andthen m
& = &
m \conjun (m^{\perp} \disjun k) \\
& \leq &
m \conjun (m^{\perp} \disjun (m\sasaki n)) \\
& = &
m \conjun (m^{\perp} \disjun (m^{\perp} \disjun (m\conjun n))) \\
& = &
m \conjun (m^{\perp} \disjun (m\conjun n)) \\
& = &
m \conjun n \qquad 
   \mbox{by orthomodularity, since }m \leq m\conjun n \\
& \leq &
n.
\end{array}$$
}

%% ``op-reindexing'' does not preserve the Sasaki hook. For a
%% counterexample, consider the kernels $\kappa,\Delta\colon X
%% \rightarrowtail X \oplus X$ in a dagger category with dagger
%% biproducts and kernels (such as $\Hilb$), and the cokernel $\pi:X
%% \oplus X \twoheadrightarrow X$. Then $\kappa \sasaki \Delta =
%% \kappa^\perp \disjun (\kappa \conjun \Delta) = \kappa' \disjun 0 = \kappa'$,
%% so that $\exists_{\pi}(\kappa \sasaki \Delta) = \Im(\pi \after
%% \kappa') = \Im(0)=0$. But $\exists_{\pi}(\kappa) = \Im(\pi \after
%% \kappa) = \Im(\idmap) = \idmap$ and $\exists_{\pi}(\Delta) = \Im(\pi
%% \after \Delta) = \Im(\idmap) = \idmap$, so that $\exists_{\pi}(\kappa)
%% \sasaki \exists_{\pi}(\Delta) = \idmap \sasaki \idmap = \idmap^\perp
%% \disjun (\idmap \conjun \idmap) = 0 \disjun \idmap = \idmap$. Hence
%% $\exists_{\pi}(\kappa \sasaki \Delta) \neq \exists_{\pi}(\kappa)
%% \sasaki \exists_{\pi}(\Delta)$ (unless $X \cong 0$).

%% The same counterexample shows that ``op-reindexing'' does not preserve
%% the `and-then' operation either. For $\kappa \andthen \Delta = \Delta
%% \conjun (\Delta^\perp \disjun \kappa) = \Delta \conjun 0 = 0$, so that
%% $\exists_{\pi}(\kappa \andthen \Delta) = \Im(\pi \after 0) = 0$.  But
%% $\exists_{\pi}(\kappa) \andthen \exists_{\pi}(\Delta) = \idmap
%% \andthen \idmap = \idmap \conjun (\idmap^\perp \disjun \idmap) = \idmap
%% \conjun \idmap = \idmap$. Hence $\exists_{\pi}(\kappa \andthen \Delta)
%% \neq \exists_{\pi}(\kappa) \andthen \exists_{\pi}(\Delta)$ (unless $X
%% \cong 0$).

As we have seen, substitution functors $f^{-1}$ in dagger kernel
categories have left adjoints $\exists_f$. It is natural to ask if
they also have right adjoints $\forall_f$. The next result says that
existence of such adjoints $\forall_f$ makes the logic Boolean. 

\begin{proposition}
\label{ForallProp}
  Suppose there are right adjoints
  $\forall_f$ to $f^{-1} \colon \KSub(Y) \to \KSub(X)$ for each $f \colon X \to
  Y$ in a dagger kernel category. Then each $\KSub(X)$ is a Boolean algebra.
\end{proposition}

\begin{proof}
  \cite[Lemma~A1.4.13]{Johnstone02} For $k,l \in \KSub(X)$, define
  implication $(k \Rightarrow l) = \forall_k(k^{-1}(l)) \in
  \KSub(X)$. Then for any $m \in \KSub(X)$:
$$\begin{array}{rcl}
m \;\leq\; \forall_k(k^{-1}(l)) = (k \Rightarrow l)
& \Longleftrightarrow & 
k^{-1}(m) \;\leq\; k^{-1}(l) \\
& \Longleftrightarrow & 
m \conjun k = k \after k^{-1}(m) \;\leq\; l,
\end{array}$$

\noindent where the last equivalence holds because $k\after -$ is left
adjoint to $k^{-1}$, since $k$ is a kernel. Hence $\KSub(X)$ is a
Heyting algebra, and therefore distributive.  By
Proposition~\ref{OrthmodularityProp} we know that it is also
orthomodular. Hence each $\KSub(X)$ is a Boolean algebra.  
\end{proof}

These universal quantifiers $\forall_f$ do not exist in general
because not all kernel posets $\KSub(X)$ are Boolean algebras.  For a
concrete non-example, consider the lattice $\KSub(\field{C}^2)$ in the
category $\Hilb$---where $\field{C}$ denotes the complex numbers. 
Consider the kernel subobjects represented by
$$\kappa_1 \colon \field{C} \to \field{C}^2, 
\qquad
\kappa_2 = (\kappa_1)^\perp \colon \field{C} \to \field{C}^2,
\qquad
\Delta = \tuple{\idmap}{\idmap} \colon \field{C} \to \field{C}^2.$$

\noindent Since we can write each $(z,w)\in\field{C}^{2}$ as $(z,w) =
\Delta(z,z) + \kappa_{2}(w-z)$ we get $\Delta \disjun \kappa_2 = 1$ in
$\KSub(\field{C}^2)$.  This yields a counterexample to
distributivity:
$$\kappa_1 \conjun (\Delta \disjun \kappa_2) 
= 
\kappa_1 \conjun 1 
=
\kappa_1 
\neq 
0 
= 
0 \disjun 0 
= 
(\kappa_1 \conjun \Delta) \disjun (\kappa_1 \conjun \kappa_2).$$

We now turn to a more systematic study of Booleanness. As we have
seen, the categories \Rel, \PInj and $\widehat{B}$ (for a Boolean
algebra $B$) are Boolean, but \Hilb and \PHilb are not. The following
justifies the name ``Boolean''.

\begin{theorem}
\label{BooleanLem}
A dagger kernel category is Boolean if and only if each orthomodular
lattice $\KSub(X)$ is a Boolean algebra.
\end{theorem}

\begin{proof}
We already know that each poset $\KSub(X)$ is an orthomodular lattice,
with bottom $0$, top $1$, orthocomplement $(-)^{\perp}$ (by
Lemma~\ref{KerLem}), intersections $\conjun$ (by
Lemma~\ref{KernelCompLem}), and joins $m\disjun n = (m^{\perp} \conjun
n^{\perp})^{\perp}$. What is missing is distributivity $m \conjun (n
\disjun k) = (m \disjun n) \conjun (m \disjun k)$. We show that the
latter is equivalent to the Booleanness requirement $m \conjun n = 0
\Rightarrow m \orthogonal n$. Recall: $m\orthogonal n$ iff $n^{\dag}
\after m = 0$ iff $m\leq n^{\perp} = \ker(n^{\dag})$.

First, assume Booleanness. In any lattice one has
$m \wedge (n \vee k) \geq (m \wedge n) \vee (m \wedge k)$.
For the other inequality, notice that
$$(m \conjun (m \conjun n)^{\perp}) \conjun n
=
(m \conjun n) \conjun (m \conjun n)^{\perp}
=
0.$$

\noindent Hence $m \conjun (m \conjun n)^{\perp} \leq n^{\perp}$. Similarly,
$m \conjun (m \conjun k)^{\perp} \leq k^{\perp}$. So
$$m \conjun (m \conjun n)^{\perp} \conjun (m \conjun k)^{\perp} 
\leq 
n^{\perp} \conjun k^{\perp}
=
(n \disjun k)^{\perp},$$

\noindent and therefore
$$m \conjun (m \conjun n)^{\perp} \conjun (m \conjun k)^{\perp} 
   \conjun (n\disjun k)
=
0.$$

\noindent But then we are done by using Booleanness again:
$$m \conjun (n\disjun k)
\leq
((m \conjun n)^{\perp} \conjun (m \conjun k)^{\perp})^{\perp}
=
(m \conjun n) \disjun (m \conjun k).$$

The other direction is easier: if $m\conjun n = 0$, then
$$\begin{array}{rcl}
m
\hspace*{\arraycolsep}=\hspace*{\arraycolsep}
m \conjun 1
& = &
m \conjun (n\disjun n^{\perp}) \\
& = &
(m \conjun n) \disjun (m \conjun n^{\perp}) 
   \quad\mbox{by distributivity} \\
& = &
0 \disjun (m \conjun n^{\perp}) 
\hspace*{\arraycolsep}=\hspace*{\arraycolsep}
m \conjun n^{\perp},
\end{array}$$

\noindent whence $m \leq n^{\perp}$. 

\end{proof}

%% Using the previous theorem we can strengthen
%% Proposition~\ref{ForallProp}, as follows.

%% \begin{corollary}
%%   A dagger kernel category that is not Boolean cannot have right
%%   adjoints $\forall_f$ to $f^{-1}$ for all morphisms $f$.
%% \end{corollary}
%% \begin{proof}
%%   If the category is not Boolean, by the previous theorem there is an
%%   object $X$ such that $\KSub(X)$ is an orthomodular lattice that is
%%   not Boolean. That means that $\KSub(X)$ is not distributive, 
%%   and hence cannot be a Heyting algebra. Hence the existence of
%%   universal quantifiers would contradict Proposition~\ref{ForallProp}.
%%   
%% \end{proof}

The Booleanness property can be strengthened in the
following way.

\begin{proposition}
\label{BooleanStrengthenedProp}
The Booleanness requirement $m\conjun n = 0 \Rightarrow m\leq
n^{\perp}$, for all kernels $m,n$, is equivalent to the following:
for each pullback of kernels:
$$\begin{array}{rcrcl}
\raisebox{1.5em}{$\xymatrix{
P\ar[r]^-{p}\ar[d]_{q}\pullback & N\ar@{ |>->}[d]^{n} \\
M\ar@{ |>->}[r]_-{m} & X
}$}
& \qquad\mbox{one has}\qquad &
n^{\dag} \after m
& = &
p \after q^{\dag}.
\end{array}$$
\end{proposition}

\begin{proof}
It is easy to see that the definition of Booleanness is the special case
$P=0$. For the converse, we put another pullback on top of the
one in the statement:
$$\xymatrix{
0\ar[r]\ar[d]\pullback & P^{\perp}\ar@{ |>->}[d]^{p^{\perp}} \\
P\ar@{ |>->}[r]^-{p}\ar@{ |>->}[d]_{q}\pullback & N\ar@{ |>->}[d]^{n} \\
M\ar@{ |>->}[r]_-{m} & X
}$$

\noindent We use that $p,q$ are kernels by Lemma~\ref{PullbackLem}. We
see $m \conjun (n \after p^{\perp}) = 0$, so by Booleanness
we obtain:
$$\begin{array}{rcl}
m \leq (n \after p^{\perp})^{\perp} 
& = &
\ker\Big((n \after \ker(p^{\dag}))^{\dag}\Big) \\
& = &
\ker(\coker(p) \after n^{\dag}) \\
& = &
(n^{\dag})^{-1}(p),
\end{array}$$

\noindent where the pullback is as described in
Lemma~\ref{PullbackLem}. Hence there is a map $\varphi\colon M
\rightarrow P$ with $p \after \varphi = n^{\dag} \after m$. This means
that $\varphi = p^{\dag} \after p \after \varphi = p^{\dag} \after n^{\dag}
\after m = (n \after p)^{\dag} \after m = (m \after q)^{\dag} \after m
= q^{\dag} \after m^{\dag} \after m = q^{\dag}$. Hence we have
obtained $p \after q^{\dag} = n^{\dag} \after m$, as
required. 

\end{proof}

\begin{definition}
\label{KcKDef}
Let \Cat{D} be a Boolean dagger kernel category. We write
$\Cat{D}_{kck}$ for the category with the same objects as \Cat{D};
morphisms $X\rightarrow Y$ in $\Cat{D}_{kck}$ are cokernel-kernel
pairs $(c, k)$ of the form $\xyinline{X \ar@{-|>}^-{c}[r] & \,\bullet\,
\ar@{ |>->}^-{k}[r] & Y}$. The identity $X\rightarrow X$ is
$\xyinline{X \ar@{-|>}^-{\idmap[]}[r] & X \ar@{ |>->}^-{\idmap[]}[r] & X}$, and
composition of $\xyinline{X \ar@{-|>}^-{c}[r] & M \ar@{ |>->}^-{k}[r]
& Y}$ and $\xyinline{Y \ar@{-|>}^-{d}[r] & N \ar@{ |>->}^-{l}[r] & Z}$
is the pair $(q^{\dag} \after c,
l \after p)$ obtained via the pullback:
\begin{equation}
\label{eq:kckcomposition}
\vcenter{\xymatrix{
& P\ar@{ |>->}[r]^-{p}\ar@{ |>->}[d]_{q}\pullback & 
   N\ar@{ |>->}[d]^{d^{\dag}}\ar@{ |>->}[r]^-{l} & Z \\
X\ar@{-|>}[r]_-{c} & M\ar@{ |>->}[r]_-{k} & Y
}}
\end{equation}

\noindent To be precise, we identity $(c,k)$ with $(\varphi\after c,
k\after \varphi^{-1})$, for isomorphisms $\varphi$. 
\end{definition}

The reader may have noticed that this construction generalises
the definition of \PInj. Indeed, now we can say $\PInj = \Rel_{kck}$.

\begin{theorem}
\label{KcKThm}
The category $\Cat{D}_{kck}$ as described in Definition~\ref{KcKDef}
is again a Boolean dagger kernel category, with a functor
$D\colon\Cat{D}_{kck}\rightarrow\Cat{D}$ that is a morphism of \DCK,
and a change-of-base situation (pullback):
$$\xymatrix@R-.75pc{
\KSub(\Cat{D}_{kck})\ar[d]\ar[rr] & & \KSub(\Cat{D})\ar[d] \\
\Cat{D}_{kck}\ar[rr]^-{D} & & \Cat{D}
}$$

\noindent Moreover, in $\Cat{D}_{kck}$ one has:
$$\mbox{kernel}
\;=\;
\mbox{dagger mono}
\;=\;
\mbox{mono}
\;=\;
\mbox{zero-mono},$$

\noindent and $\Cat{D}_{kck}$ is universal among such categories.
\end{theorem}

\begin{proof}
The obvious definition $(c,k)^{\dag} = (k^{\dag}, c^{\dag})$ yields an
involution on $\Cat{D}_{kck}$. The zero object $0\in\Cat{D}$ is also a
zero object $0\in\Cat{D}_{kck}$ with zero map $\xyinline[@C-2ex]{X
  \ar@{-|>}[r] & \;0\; \ar@{ |>->}[r] & Y}$ consisting of a
cokernel-kernel pair. A map $(c,k)$ is a dagger mono if and only if
$(c,k)^{\dag} \after (c,k) = (k^{\dag},k)$ is the identity; this means
that $k=\idmap$.

The kernel of a map $(d,l) = (\xyinline[@C-2ex]{Y \ar@{-|>}^-{d}[r] &
  \;N\; \ar@{ |>->}^-{l}[r] & Z})$ in $\Cat{D}_{kck}$ is $\ker(d,l) =
(\xyinline[@C-2ex]{N^{\perp} \ar@{-|>}^-{\idmap}[r] & \;N^{\perp}\; \ar@{
    |>->}^-{(d^{\dag})^{\perp}}[r] & Y})$, so that $\ker(d,l)$ is a
dagger mono and $(d,l) \after \ker(d,l) = 0$. If also $(d,l) \after
(c,k) = 0$, then $k \conjun d^{\dag} = 0$ so that by Booleanness, $k
\leq (d^{\dag})^{\perp}$, say via $\varphi \colon M\rightarrow
N^{\perp}$ with $(d^{\dag})^{\perp} \after \varphi = k$. Then we
obtain a mediating map $(c,\varphi) = (\xyinline{X \ar@{-|>}^-{c}[r] &
  \;M\; \ar@{ |>->}^-{\varphi}[r] & N^{\perp}})$ which satisfies
$\ker(d,l) \after (c,\varphi) = (\idmap, (d^{\dag})^{\perp}) \after
(c,\varphi) = (c, (d^{\dag})^{\perp} \after \varphi) = (c,k)$.
It is not hard to see that maps of the form $(\idmap, m)$ in
$\Cat{D}_{kck}$ are kernels, namely of the cokernel $(m^{\perp},\idmap)$.

The intersection of two kernels $(\idmap,m) = (\xyinline[@C-2ex]{M
  \ar@{=}[r] & \;M\; \ar@{ |>->}^-{m}[r] & X})$ and $(\idmap,n) =
(\xyinline[@C-2ex]{N \ar@{=}[r] & \;N\; \ar@{ |>->}^-{n}[r] & X})$ in
$\Cat{D}_{kck}$ is the intersection $m\conjun n \colon P
\rightarrowtail X$ in \Cat{D}, with projections $(\xyinline[@C-2ex]{P
  \ar@{=}[r] & \;P\; \ar@{ |>->}[r] & M})$ and $(\xyinline[@C-2ex]{P
  \ar@{=}[r] & \;P\; \ar@{ |>->}[r] & N})$. Hence if the intersection
of $(\idmap,m)$ and $(\idmap,n)$ in $\Cat{D}_{kck}$ is 0, then so is
the intersection of $m$ and $n$ in \Cat{D}, which yields $n^{\dag}
\after m = 0$. But then in $\Cat{D}_{kck}$, $(\idmap,n)^{\dag} \after
(\idmap,m) = (n^{\dag},\idmap) \after (\idmap,m) = 0$. Hence
$\Cat{D}_{kck}$ is also Boolean.

Finally, there is a functor $\Cat{D}_{kck} \rightarrow \Cat{D}$
by $X\mapsto X$ and $(c,k) \mapsto k\after c$. Composition is preserved
by Proposition~\ref{BooleanStrengthenedProp}, since for maps as
in Definition~\ref{KcKDef},
$$\begin{array}{rcl}
(d,l) \after (c,k)
\hspace*{\arraycolsep}=\hspace*{\arraycolsep}
(q^{\dag} \after c, l \after p)
& \longmapsto &
l \after p \after q^{\dag} \after c 
\hspace*{\arraycolsep}=\hspace*{\arraycolsep}
(l \after d) \after (k \after c).
\end{array}$$

\noindent We have already seen that $\KSub(X)$ in $\Cat{D}_{kck}$
is isomorphic to $\KSub(X)$ in \Cat{D}. This yields the change-of-base
situation.

We have already seen that kernels and dagger monos coincide. We now
show that they also coincide with zero-monos. So let $(d,l)\colon Y
\rightarrow Z$ be a zero-mono. This means that $(d,l) \after (c,k) = 0
\Rightarrow (c,k) = 0$, for each map $(c,k)$. Using
diagram~(\ref{eq:kckcomposition}), this means: $d^{\dag} \conjun k = 0
\Rightarrow k=0$. By Booleanness, the antecedent $d^{\dag} \conjun k =
0$ is equivalent to $k \leq (d^{\dag})^{\perp} = \ker(d)$, which means
$d\after k = 0$. Hence we see that $d$ is zero-monic in \Cat{D}, and
thus an isomorphism (because it is already a cokernel).

Finally, let \Cat{E} be a Boolean dagger kernel category in
which zero-monos are kernels, with a functor $F\colon \Cat{E}
\rightarrow \Cat{D}$ in \DCK. Every morphism $f$ in \Cat{E} factors as
$f=i_f \after e_f$ for a kernel $i_f$ and a cokernel $e_f$. Hence
$G\colon \Cat{E} \to \Cat{D}_{kck}$ defined by $G(X) = F(X)$ and $G(f) =
(e_f,i_f)$ is the unique functor satisfying $F = D \after G$. 
\end{proof}

\auxproof{
We check that $\dag$ preserves composition:
$$\begin{array}{rcll}
\big((d,l) \after (c,k)\big)^{\dag}
& = &
(q^{\dag} \after c, l \after p)^{\dag}
   & \mbox{see the pullback in the proposition} \\
& = &
\big((l \after p)^{\dag}, (q^{\dag} \after c)^{\dag}\big) \\
& = &
(p^{\dag} \after l^{\dag}, c^{\dag} \after q) \\
& = &
(k^{\dag}, c^{\dag}) \after (l^{\dag}, d^{\dag})
   & \mbox{see the diagram below} \\
& = &
(c,k)^{\dag} \after (d,l)^{\dag}.
\end{array}$$

$$\xymatrix{
& P\ar@{ |>->}[r]^-{q}\ar@{ |>->}[d]_{p}\pullback & 
   M\ar@{ |>->}[d]^{k}\ar@{ |>->}[r]^-{c^{\dag}} & X \\
Z\ar@{-|>}[r]^-{l^{\dag}} & N\ar@{ |>->}[r]^-{d^{\dag}} & Y
}$$
}

\section{Ordering homsets}\label{HomsetOrderSec}

This section shows that homsets in dagger kernel categories
automatically carry a partial order. However, this does not make the
categories order enriched, because the order is not preserved by all
morphisms.

\begin{definition} 
\label{OrderDef}
Let $f,g\colon X \rightarrow Y$ be parallel morphisms in a dagger
kernel category. After factorising them as $f=i_f \after m_f \after
(i_{f^\dag})^\dag$ and $g=i_g \after m_g \after (i_{g^\dag})^\dag$
like in~(\ref{ImageCoimageEqn}) we can define $f \leq g$ if and only
if there are (necessarily unique, dagger monic) $\varphi\colon \Im(f)
\to \Im(g)$ and $\psi \colon \Im(f^{\dag}) \to \Im(g^{\dag})$, so that
in the diagram
  \begin{equation}
  \label{eq:orderdiagram}
   \vcenter{\xymatrix@C+2ex@R-3ex{
    &   \Im(f^{\dag}) \ar^-{m_f}[r]
    &   \Im(f)\ar@{ |>->}^-{i_f}[dr] \ar@{-->}_-{\varphi}[dd]
    &   \\
        X \ar@{-|>}^-{(i_{f^\dag})^{\dag}}[ur] 
          \ar@{-|>}_-{(i_{g^\dag})^{\dag}}[dr]
    &&& Y \\
    &   \Im(g^{\dag}) \ar_-{m_g}[r] \ar@{-->}_-{\psi^{\dag}}[uu]
    &   \Im(g) \ar@{ |>->}_-{i_g}[ur]
  }}
  \end{equation}
  one has
  \[
    \psi^{\dag} \after (i_{g^\dag})^{\dag} = (i_{f^\dag})^{\dag}
    \quad
    \varphi\after m_f = m_g \after \psi
    \quad
    \varphi^{\dag} \after m_g = m_f \after \psi^{\dag}
    \quad
    i_g \after \varphi = i_f.
  \]
\end{definition}

\begin{lemma}
  The relation $\leq$ is a partial order on each homset of a dagger
  kernel category, with the zero morphism as least element. 
\end{lemma}

\begin{proof}
  Reflexivity is easily established by taking $\varphi=\idmap$ and
  $\psi=\idmap$ in~(\ref{eq:orderdiagram}).  For transitivity, suppose
  that $f \leq g$ via $\varphi$ and $\psi$, and that $g \leq h$ via
  $\alpha$ and $\beta$. Then the four conditions in the previous
  definition are fulfilled by $\alpha \after \varphi$ and $\psi \after
  \beta$, so that $f \leq h$.  Finally, for anti-symmetry, suppose
  that $f \leq g$ via $\varphi$ and $\psi$, and that $g \leq f$ via
  $\alpha$ and $\beta$.  Then $i_{f} \after \alpha \after \varphi =
  i_{g} \after \varphi = i_{f}$, so that $\alpha \after \varphi =
  \idmap$. Similarly, $\beta \after \psi = \idmap$. By
  Lemma~\ref{KerLem}, $\alpha$ is a dagger mono so that $\alpha^{\dag}
  = \alpha^{\dag} \after \alpha \after \varphi = \varphi$. Similarly,
  $\beta^{\dag} = \psi$, and thus:
$$\begin{array}{rcl}
f
\hspace*{\arraycolsep}=\hspace*{\arraycolsep}
i_{f} \after m_{f} \after (i_{f^\dag})^{\dag}
\hspace*{\arraycolsep}=\hspace*{\arraycolsep}
i_{f} \after \alpha \after \varphi \after m_{f} \after (i_{f^\dag})^{\dag}
& = &
i_{g} \after m_{g} \after \psi \after (i_{f^\dag})^{\dag} \\
& = &
i_{g} \after m_{g} \after \beta^{\dag} \after (i_{f^\dag})^{\dag} \\
& = &
i_{g} \after m_{g} \after (i_{g^\dag})^{\dag} \\
& = &
g.
\end{array}$$

\noindent Finally, for any $f$ we have $0 \leq f$ by taking $\varphi =
\psi = 0$ in~(\ref{eq:orderdiagram}).  
\end{proof}

\begin{lemma}
\label{OrderPreservationLem}
If $f\leq g$, then:
\begin{enumerate}
\item $(k\after f) \leq (k\after g)$ for a kernel $k$;

\item $(f\after c) \leq (g\after c)$ for a cokernel $c$;

\item $f^{\dag} \leq g^{\dag}$.
\end{enumerate}
\end{lemma}

\begin{proof}
  The first two points are obvious. The third one then follows because
  $(m_{f})^{\dag} = m_{f^{\dag}}$ as shown in 
  Section~\ref{ImageCoimageSec}. 
\end{proof}

\begin{example}
\label{OrderEx}
We describe the situation in \PInj, \Rel and \Hilb, using the
factorisations from Example~\ref{ImageCoimageEx}.

Two parallel maps $\smash{f = (X \stackrel{f_1}{\leftarrowtail} F
  \stackrel{f_2}{\rightarrowtail} Y)}$ and $\smash{g = (X
  \stackrel{g_1}{\leftarrowtail} G \stackrel{g_2}{\rightarrowtail}
  Y)}$ in \PInj satisfy $f\leq g$ if and only if there are
$\varphi,\psi\colon F\rightarrow G$ in:
$$\xymatrix@C+2ex@R-4ex{
& F\ar@{=}[r] & F\ar@{ |>->}[dr]^-{f_2}\ar@{-->}[dd]_{\varphi} & \\
X\ar@{-|>}[ur]^-{(f_{1})^{\dag}}\ar@{-|>}[dr]_{(g_{1})^{\dag}} & & & Y \\
& G\ar@{=}[r]\ar@{-->}[uu]_{\psi^{\dag}} & G\ar@{ |>->}[ur]_{g_2}
}$$

\noindent This means $\varphi=\psi$ and $g_{i} \after \varphi =
f_{i}$, for $i=1,2$, so that we obtain the usual order (of one partial
injection extending another).

Next, $R\leq S$ for $\smash{R = (X \stackrel{r_1}{\leftarrow} R
  \stackrel{r_2}{\rightarrow} Y)}$ and $\smash{S = (X
  \stackrel{s_1}{\leftarrow} S \stackrel{s_2}{\rightarrow}
  Y)}$ in \Rel means:
$$\xymatrix@C+2ex@R-4ex{
& \Im(r_{1})\ar[r]^-{R} & \Im(r_{2})\ar@{ |>->}[dr]\ar@{-->}[dd]_{\varphi} & \\
X\ar@{-|>}[ur]\ar@{-|>}[dr] & & & Y \\
& \Im(s_{1})\ar[r]_-{S}\ar@{-->}[uu]_{\psi^{\dag}} & \Im(s_{1})\ar@{ |>->}[ur]
}$$

\noindent Commutation of the triangles means $\Im(r_{1})\subseteq\Im(s_{1})$
and $\Im(r_{2})\subseteq\Im(s_{2})$. The equations for the square in the
middle say that:
$$R(x,y) \Leftrightarrow S(x,y)
\quad\mbox{for all}\quad
\left\{\begin{array}{l}
(x,y)\in\Im(r_{1})\times\Im(s_{2}) \\
(x,y)\in\Im(r_{2})\times\Im(s_{1}).
\end{array}\right.$$

\noindent This means $R\subseteq S$, as one would expect.

The order on the homsets of the category $\Hilb$ can be characterized
as follows \cite[Example~5.1.10]{Heunen10a}:
$f \leq g$ for $f,g \colon X \to Y$ if and only if
$g = f+f'$ for some $f' \colon X \to Y$ with $\Im(f)$ and
$\Im(f^\dag)$ orthogonal to $\Im(f')$ and $\Im((f')^\dag)$,
respectively. 
To see this, suppose that $g=f+f'$ as above. Then $\Im(g)$ is the
direct sum of $\Im(f)$ and $\Im(f')$, and likewise $\Im(g^\dag) =
\Im(f^\dag) \oplus \Im((f')^\dag$. Moreover, $m_g$ is the
direct sum of $m_f$ and $m_{f'}$. Therefore, taking $\psi =
\varphi = \kappa_1$ makes diagram~\eqref{eq:orderdiagram} commute, so
that $f \leq g$. 
Conversely, suppose that $f \leq g$, so that
diagram~\eqref{eq:orderdiagram} commutes. Then the cotuple $[\varphi,
\varphi^\perp] \colon \Im(f) \oplus \Im(f)^\perp \to \Im(g)$ is an
isomorphism, and so is the cotuple $[\psi, \psi^\perp]$.
Since $\varphi^\dag \after m_g = m_f \after \psi^\dag$, there is a
morphism $n$ making the following diagram commute: 
\[\xymatrix@R-2ex{
    \Im(f^\dag)^\perp \ar@{ |>->}_-{\ker(\psi^\dag)=\psi^\perp}[d]
    \ar@{-->}^-{n}[rr] 
 && \Im(f)^\perp \ar@{ |>->}^{\varphi^\perp = \ker(\varphi^\dag)}[d] \\
    \Im(g^\dag) \ar^-{m_g}[rr] \ar@{-|>}_-{\psi^\dag}[d]
 && \Im(g) \ar@{-|>}^{\varphi^\dag}[d] \\
    \Im(f^\dag) \ar_-{m_f}[rr] 
 && \Im(f).
}\]
Now, taking
\[\xymatrix{
   f' = \Big( X \ar^-{(i_{g^\dag})^\dag}[r]
   & \Im(g^\dag) \ar^-{(\psi^\perp)^\dag}[r]
   & \Im(f^\dag)^\perp \ar^-{n}[r]
   & \Im(f)^\perp \ar^-{\varphi^\perp}[r]
   & \Im(g) \ar^-{i_g}[r]
   & Y \Big)
}\]
fulfills $g=f+f'$, and $\Im(f)$ and $\Im(f^\dag)$ are orthogonal to
$\Im(f')$ and $\Im((f')^\dag)$, respectively. 
\end{example}

In Hilbert spaces there is a standard correspondence between
self-adjoint idempotents and closed subsets. Recall that an endomap
$p\colon X\rightarrow X$ is self-adjoint if $p^{\dag} = p$ and
idempotent if $p \after p = p$. In the current, more general, setting
this works as follows, using the order on homsets.

\begin{proposition}
\label{KernelEndomapProp}
The ``effect''\footnote{The name ``effect'' was chosen because of
connections to effect algebras~\cite{DvurecenskijP00}. For example, 
in the so-called standard effect algebra of a Hilbert
space~\cite{FoulisGB94}, an effect corresponds a positive operator
beneath the identity.} mapping 
$m\mapsto \effect{m} 
\,\smash{\stackrel{\textrm{def}}{=}}\, m \after m^{\dag}$ from
Section~\ref{LogicSec} yields an order isomorphism:
$$\begin{array}{rcl}
\KSub(X)
& \cong &
\set{p\colon X\rightarrow X}{p^{\dag} = p\leq \idmap} \\
& \cong &
\set{p\colon X\rightarrow X}{p^{\dag} = p\after p = p\leq \idmap} \\
& \smash{\stackrel{(*)}{\cong}} &
\set{p\colon X\rightarrow X}{p^{\dag} = p\after p = p},
\end{array}$$

\noindent where the marked isomorphism holds if zero-epis are epis
(like in \Hilb).
\end{proposition}
\begin{proof}
  Clearly, $\effect{m} = m \after m^{\dag}$ is a self-adjoint idempotent. 
  It satisfies $\effect{m} \leq \idmap$ via:
$$\xymatrix@C+2ex@R-4ex{
& M\ar@{=}[r] & M\ar@{ |>->}[dr]^-{m}\ar@{-->}[dd]_{m} & \\
X\ar@{-|>}[ur]^-{m^{\dag}}\ar@{=}[dr] & & & X \\
& X\ar@{=}[r]\ar@{-->}[uu]_{m^{\dag}} & X\ar@{=}[ur] 
}$$

\noindent where the kernel $m\colon M\rightarrowtail X$ is a dagger mono
so that $\Im(\effect{m}) = M$.

This mapping $\effect{-}\colon \KSub(X) \rightarrow \set{p}{p^{\dag} =
  p \leq \idmap}$ is surjective: if $p\colon X\rightarrow X$ is a
self-adjoint with $p\leq\idmap$ then we first note that the
factorisation from~(\ref{ImageCoimageEqn}) yields $p = i_{p} \after
m_{p} \after (i_{p})^{\dag}$. By Definition~\ref{OrderDef} there are
$\varphi,\psi\colon \Im(p)\rightarrow X$ with $\psi^{\dag} =
(i_{p})^{\dag}$, $\varphi \after m_{p} = \psi$, $\varphi^{\dag} =
m_{p} \after \psi^{\dag}$ and $\varphi = i_{p}$. This yields
$\psi=i_{p}$ and $m_{p}=\idmap$.  Hence $p = i_{p} \after
(i_{p})^{\dag} = \effect{i_{p}}$, so that $p$ is automatically
idempotent. This establishes the second isomorphism.

The mapping $\effect{-}$ preserves and reflects the order. If $m\leq n$
in $\KSub(X)$, say via $\varphi\colon M\rightarrow N$ with $n\after\varphi
= m$, then $\effect{m} \leq \effect{n}$ via:
$$\xymatrix@C+2ex@R-4ex{
& M\ar@{=}[r] & M\ar@{ |>->}[dr]^-{m}\ar@{-->}[dd]_{\varphi} & \\
X\ar@{-|>}[ur]^-{m^{\dag}}\ar@{-|>}[dr]_{n^{\dag}} & & & X \\
& N\ar@{=}[r]\ar@{-->}[uu]_{\varphi^{\dag}} & N\ar@{ |>->}[ur]_{n}
}$$

\noindent Conversely, if $\effect{m} \leq \effect{n}$, say via $\varphi\colon
M\rightarrow N$ and $\psi\colon M\rightarrow N$, then $n\after \varphi = m$
so that $m\leq n$ in $\KSub(X)$.

Finally, if zero-epis are epis, we write for a self-adjoint idempotent
$p$,
$$i_{p} \after e_{p} 
= 
p 
= 
p \after p 
= 
p^{\dag} \after p 
=
(e_{p})^{\dag} \after (i_{p})^{\dag} \after i_{p} \after e_{p} 
=
(e_{p})^{\dag} \after e_{p},$$ 

\noindent and obtain $i_{p} = (e_{p})^{\dag}$. Hence $p = \effect{i_{p}}$
and thus $p\leq \idmap$. 
\end{proof}

\section{Completeness and atomicity of kernel posets}
\label{ComplAtomSec} 

In traditional quantum logic, orthomodular lattices are usually
considered with additional properties, such as completeness and
atomicity~\cite{Piron76}. This section considers how these
requirements on the lattices $\KSub(X)$ translate to categorical
properties. For convenience, let us recall the following standard
order-theoretical definitions.completeness

\begin{definition}
  For elements $x,y$ of a poset, we say that $y$ covers $x$ when
  $x<y$ and $x \leq z < y$ implies $z=x$ (where $z<y$ if and only if
  $z \leq y$ and $z \neq y$). An element $a$ of a poset with
  least element 0 is called an \emph{atom} when it covers
  0. Equivalently, an atom cannot be expressed as a join of strictly
  smaller elements. Consequently, $0$ is not an atom. A poset is
  called \emph{atomic} if for any $x \neq 0$ in it there exists an
  atom $a$ with $a \leq x$. Finally, a lattice is \emph{atomistic}
  when every element is a join of atoms~\cite{DaveyP90}.
\end{definition}

\begin{proposition}
\label{prop:atoms}
  For an arbitrary object $I$ in a dagger kernel category, the
  following are equivalent: 
  \begin{enumerate}
  \item $\idmap[I]=1$ is an atom in $\KSub(I)$;
  \item $\KSub(I)=\{0,1\}$;
  \item each nonzero kernel $x\colon I \rightarrowtail X$ is an atom
    in $\KSub(X)$.
  \end{enumerate}
\end{proposition}

\begin{proof}
  For the implication $(1) \Rightarrow (2)$, let $m$ be a kernel into
  $I$. Because $m \leq \idmap[I]$ and the latter is an atom, we have
  that $m=0$ or $m$ is isomorphism. Thus $\KSub(I)=\{0,1\}$.

  To prove $(2) \Rightarrow (3)$, suppose that $m \leq x$ for kernels
  $m\colon M \rightarrowtail X$ and $x\colon I \rightarrowtail X$.
  Say $m = x \after \varphi$ for $\varphi\colon M \rightarrowtail I$.
  Then $\varphi$ is a kernel by Lemma~\ref{KerLem}.  Since
  $\KSub(I)=\{0,1\}$, either $\varphi$ is zero or $\varphi$ is
  isomorphism. Hence either $m=0$ or $m=x$ as subobjects.  So $x$ is
  an atom.  Finally, $(3) \Rightarrow (1)$ is trivial.  
\end{proof}

\begin{definition}
\label{def:KSubsimple}
  If $I$ satisfies the conditions of the previous lemma, we call it a
  \emph{$\KSub$-simple} object. (Any simple object in the usual sense
  of category theory is $\KSub$-simple.)

  Similarly, let us call $I$ a \emph{$\KSub$-generator} if
  $f=g\colon X \to Y$ whenever $f \after x = g \after x$ for all
  kernels $x\colon I \rightarrowtail X$.  
  (Any $\KSub$-generator is a generator in the usual sense of category
  theory.) 
\end{definition}

\begin{example}
\label{GeneratorEx}
The objects $1 \in \PInj$, $1 \in \Rel$, $\field{C} \in \Hilb$ and
$\field{C} \in \Cat{PHilb}$ are $\KSub$-simple $\KSub$-generators.

The two-element orthomodular lattice $2$ is a generator in the
category \Cat{OMLatGal} from~\cite{Jacobs09a}, because maps
$2\rightarrow X$ correspond to elements in $X$. But $2$ is not a
$\KSub$-generator: these maps $2\rightarrow X$ are not kernels.
\end{example}

Because $1 \in \Cat{Rel}$ is a $\KSub$-simple $\KSub$-generator,
one might expect a connection between Definition~\ref{def:KSubsimple}
and the ``kernel opclassifiers'' discussed at the end of
Section~\ref{FactorisationSec}. There is, however, no apparent such
connection. For example, the object $1$ in the category $\Cat{PInj}$ is a
$\KSub$-simple $\KSub$-generator, but not a ``kernel opclassifier''.

\begin{lemma}
\label{lem:atomic}
Suppose that a dagger kernel category $\Cat{D}$ has a $\KSub$-simple
$\KSub$-generator $I$. Then beneath any nonzero element of $\KSub(X)$
lies a nonzero element of the form $x\colon I \rightarrowtail
X$. Hence $\KSub(X)$ is atomic, and its atoms are the nonzero kernels
$x \colon I \rightarrowtail X$.
\end{lemma}

\begin{proof}
Suppose $m\colon M\rightarrowtail X$ is a nonzero kernel. Since $I$ is
a $\KSub$-generator, there must be a kernel $x\colon I\rightarrowtail
M$ with $m\after x \neq 0$. By Proposition~\ref{prop:atoms} this
$m\after x$ is an atom. It satisfies $m\after x \leq m$, so we are
done. 

%%   Suppose that there is a nonzero kernel $m \colon M \to X$ through
%%   which no nonzero kernel $x \colon I \to X$ factors. 
%%   If $f \colon I \to M$ is a kernel, then
%%   $x=m \after f$ is a kernel that factors through $m$, so $x$ must be
%%   $0$. And since $m$ is mono, in that case also $f=0$. Hence $0$ is
%%   the only kernel $I \to M$.
%%   Because $I$ is a $\KSub$-generator, it follows
%%   that $\Cat{D}(M,Y)=\{0\}$ for any object $Y$. But then $m=0$, which
%%   is a contradiction. 
%%   
\end{proof}

\begin{corollary}
\label{lem:atomistic}
  If a dagger kernel category has a $\KSub$-simple $\KSub$-generator 
  $I$, then $\KSub(X)$ is atomistic for any object $X$.
\end{corollary}
\begin{proof}
  Any atomic orthomodular lattice is atomistic~\cite{Kalmbach83}. 
  % Let $m \in \KSub(X)$. We show that $m$ is the least upper bound of
  % the set $\downset_{A}m = \set{x \leq m}{x \mbox{ atom}}$. Obviously
  % $m$ is an upper bound of $\downset_{A}m$. Suppose that $x \leq n$
  % for all $x \in \downset_{A}m$. We have to prove that $m \leq n$, or
  % equivalently, $m = n \after n^\dag \after m$. This is obvious if
  % $m=0$, so assume $m\neq 0$. Since $I$ is a $\KSub$-generator, it
  % suffices to prove $m \after y = n \after n^\dag \after m \after y$
  % for all (nonzero) kernels $y\colon I \rightarrowtail M$. Now, $x=m
  % \after y$ is nonzero and thus an atom by
  % Proposition~\ref{prop:atoms}. From $x \leq m$ we get $x\in
  % \downset_{A}m$ and thus $x \leq n$, so that $m \after y = x = n
  % \after n^\dag \after x = n \after n^\dag \after m \after y$. Thus
  % $\KSub(X)$ is atomistic.  
\end{proof}

The categorical requirement of a simple generator is quite natural in
this setting, as it is also used to prove that a certain class of
dagger kernel categories embeds into $\Hilb$~\cite{Heunen09b}.

We now turn to completeness, by showing that the existence of directed
colimits ensures that kernel subobject lattices are complete. This,
too, is a natural categorical requirement in the context of
infinite-dimensionality~\cite{Heunen08a}.  Recall that a
\emph{directed colimit} is a colimit of a directed poset, considered
as a diagram. The following result can be obtained abstractly in two
steps: directed colimits in \Cat{D} yield direct colimits in slice
categories $\Cat{D}/X$, see~\cite[Vol.~2,
Prop.~2.16.3]{Borceux94}. The reflection $\KSub(X) \hookrightarrow
\Cat{D}/X$ induced by factorisation transfers these directed colimits
to $\KSub(X)$.  However, in the proof below we give a concrete
construction.

\begin{proposition}
  If a dagger kernel category $\Cat{D}$ has directed colimits, then
  $\KSub(X)$ is a complete lattice for every $X \in \Cat{D}$.  
\end{proposition}

\begin{proof}
A lattice is complete if it has directed joins
(see~\cite[Lemma~I.4.1]{Johnstone82},
or~\cite[Lemma~2.12]{Kalmbach86}), so we shall prove that $\KSub(X)$
has such directed joins. Let $(m_{i}\colon M_{i}\rightarrowtail
X)_{i\in I}$ be a directed collection in $\KSub(X)$. For $i\leq j$ we
have $m_{i} \leq m_{j}$ and thus $m_{j} \after m_{j}^{\dag} \after
m_{i} = m_{i}$.

Let $M$ be the colimit in \Cat{D} of the domains $M_i$, say with
coprojections $c_{i}\colon M_{i}\rightarrow M$. The $(m_{i}\colon
M_{i}\rightarrowtail X)_{i\in I}$ form a cocone by assumption, so
there is a unique map $m\colon M\rightarrow X$ with $m\after c_{i} =
m_{i}$. The kernel/zero-epi factorisation~(\ref{ImageEqn}) yields:
$$\xymatrix{
m = \big(M\ar@{->>}[r]|-{\circ}^-{e} & N\ar@{ |>->}[r]^-{n} & X\Big)
}$$

\noindent We claim that $n$ is the join in $\KSub(X)$ of the $m_i$.
\begin{itemize}
\item $m_{i} \leq n$ via $e\after c_{i}\colon M_{i} \rightarrow N$
satisfying $n \after (e\after c_{i}) = m \after c_{i} = m_{i}$.

\item If $m_{i} \leq k$, then $k\after k^{\dag} \after m_{i} = m_{i}$.
Also, the maps $k_{i} = k^{\dag} \after m_{i} \colon M_{i}\rightarrow K$ form
a cocone in \Cat{D} because the $m_i$ are directed and $k$ is monic: 
if $i\leq j$, then,
$$k\after k_{j} \after m_{j}^{\dag} \after m_{i}
=
k \after k^{\dag} \after m_{j} \after m_{j}^{\dag} \after m_{i}
=
k \after k^{\dag} \after m_{i}
=
k\after k_{i}.$$

\noindent As a result there is a unique $\ell\colon M\rightarrow K$ with
$\ell \after c_{i} = k_{i}$. Then $k \after \ell = m$ by uniqueness since:
$$k\after \ell \after c_{i}
=
k\after k_{i}
=
k \after k^{\dag} \after m_{i}
=
m_{i}
=
m \after c_{i}.$$

\noindent Hence we obtain $n\leq k$ by diagonal-fill-in from
Lemma~\ref{DiagonalFillInLem} in:
$$\vcenter{\xymatrix{
M\ar@{->>}[r]|-{o}^-{e}\ar[d]_{\ell} & 
   N\ar@{ |>->}[d]^{n}\ar@{..>}[dl] \\
K\ar@{ |>->}[r]_-{k} & X
}}$$
\end{itemize}
\end{proof}

\begin{example}
The categories $\PInj$, $\Rel$, $\Hilb$ and $\Cat{PHilb}$ have
directed colimits, and therefore their kernel subobject lattices are
complete orthomodular lattices.  Since they also have appropriate
generators, see Example~\ref{GeneratorEx}, each $\KSub(X)$ in $\PInj$,
$\Rel$, $\Hilb$ or $\Cat{PHilb}$ is a complete atomic atomistic
orthomodular lattice.

Any atom of a Boolean algebra $B$ is a $\KSub$-simple object in the
dagger kernel category $\widehat{B}$ from
Proposition~\ref{BAConstrProp}. But $\widehat{B}$ has a
$\KSub$-generator only if $B$ is atomistic. In that case the greatest
element 1 is a $\KSub$-generator. For if $f \after a = g \after a$ for
all $a \leq 1 \conjun x = x$ 
and $f,g \leq x \conjun y$, then, writing $\downset_{A}x = 
\setin{a}{\mathrm{Atoms}(B)}{a\leq x}$ we get:
\begin{align*}
\textstyle f
=
f \conjun x
=
f \conjun \big(\bigvee\downset_{A}x\big)
& =
\textstyle\bigvee\set{f\conjun a}{a\in\mathrm{Atoms}(B), a\leq x} \\
& =
\textstyle\bigvee\set{g\conjun a}{a\in\mathrm{Atoms}(B), a\leq x} \\
& = 
\textstyle g \conjun \big(\bigvee\downset_{A}x\big)
=
g \conjun x
= 
g.
\end{align*}
\end{example}

\section{Conclusions and future work}\label{ConclSec}

The paper shows that a ``dagger kernel category'' forms a simple but
powerful notion that not only captures many examples of interest in
quantum logic but also provides basic structure for categorical logic.
There are many avenues for extension and broadening of this work, by
including more examples (\textit{e.g.}~effect
algebras~\cite{DvurecenskijP00}) or more structure (like
tensors). Also, integrating probabilistic aspects of quantum logic is
a challenge.

\subsection*{Acknowledgements}

Thanks are due to an anonymous referee for pointing out the
reference~\cite{Crown75}. It already contains several early ideas that
are rediscovered and elaborated in the current (categorical)
setting. The category \Cat{OMLatGal} that plays a central role
in~\cite{Jacobs09a} is also already in~\cite{Crown75}.

% Several research avenues are still open: construction of dagger kernel
% categories from orthomodular lattices (like in
% Subsection~\ref{BAConstrSubsec} for Boolean algebras), or further
% investigation of the relevance of ``opfibred'' structure in this
% setting (like for ``op-classifiers'' at the end of
% Section~\ref{FactorisationSec}).

% A follow-up paper is in preparation; it extends the
% present setting with tensors (both sums $\oplus$ and products
% $\otimes$), which lead to further logical structure.

\bibliographystyle{plain}
\bibliography{dagkercat-arxiv}

\end{document}